\documentclass[12pt]{article}
\usepackage{latexsym}
\usepackage[authoryear]{natbib}

\def\R{{\rm I \mkern-2.5mu \nonscript\mkern-.5mu R}}

\usepackage{float}
\usepackage{graphicx}

\newcommand{\be}{\begin{equation}}
\newcommand{\ee}{\end{equation}}
\newcommand{\bea}{\begin{eqnarray}}
\newcommand{\eea}{\end{eqnarray}}
\newcommand{\bean}{\begin{eqnarray*}}
\newcommand{\eean}{\end{eqnarray*}}

\newcommand{\sect}[1]{\section{#1} \setcounter{equation}{0}
                      \setcounter{table}{0} \setcounter{figure}{0}}

\newtheorem{defi}{Definition}[section]

\newtheorem{lem}[defi]{Lemma}
\newtheorem{theo}[defi]{Theorem}
\newtheorem{cor}[defi]{Corollary}

\newcommand{\halmos}{\vspace{3mm} \hfill \mbox{$\Box$}}

\title{Simulation of multivariate diffusion bridges}

\author{Mogens Bladt \\
{\small Instituto de Investigacion en Matem\'aticas Aplicadas y en 
Sistemas \vspace{-1.5mm}} \\
{\small Universidad Nacional Aut\'onoma de M\'exico 
\vspace{-1.5mm}} \\ 
{\small A.P.\ 20-726\vspace{-1.5mm}} \\ 
{\small 01000 Mexico, D.F.\vspace{-1.5mm}} \\
{\small Mexico\vspace{-1.5mm}} \\
{\small bladt@sigma.iimas.unam.mx}
\and
Samuel Finch \\ 
{\small Dept.\ of Mathematical Sciences \vspace{-1.5mm}} \\
{\small University of Copenhagen\vspace{-1.5mm}} \\
{\small Universitetsparken 5\vspace{-1.5mm}} \\
{\small DK-2100 Copenhagen {\O}\vspace{-1.5mm}} \\
{\small Denmark\vspace{-1.5mm}} \\ 
{\small pjb926@math.ku.dk}
\and
Michael S\o rensen \\ 
{\small Dept.\ of Mathematical Sciences \vspace{-1.5mm}} \\
{\small University of Copenhagen\vspace{-1.5mm}} \\
{\small Universitetsparken 5\vspace{-1.5mm}} \\
{\small DK-2100 Copenhagen {\O}\vspace{-1.5mm}} \\
{\small Denmark\vspace{-1.5mm}} \\ 
{\small michael@math.ku.dk}}

\begin{document}

\maketitle

\begin{abstract}
We propose simple methods for multivariate diffusion bridge
simulation, which plays a fundamental role in simulation-based
likelihood and Bayesian inference for stochastic differential equations.
By a novel application of classical coupling methods, the new
approach generalizes a previously proposed simulation method for
one-dimensional bridges to the multi-variate setting. First a method of simulating
approximate, but often very accurate, diffusion bridges is proposed. 
These approximate bridges are used as proposal for easily
implementable MCMC algorithms that produce exact diffusion
bridges. The new method is much more generally applicable than previous methods. 
Another advantage is that the new method
works well for diffusion bridges in long intervals because
the computational complexity of the method is linear in the length of
the interval. In a simulation study the new method performs well, and its
usefulness is illustrated by an application to Bayesian estimation for
the multivariate hyperbolic diffusion model. \\ \\
{\bf Key words:} Bayesian inference; coupling; discretely sampled diffusions; likelihood 
inference; stochastic differential equation; time-reversal. \\ \\
{\bf Address for correspondance:} Michael S\o rensen, Dept.\ of
Mathematical Sciences, University of Copenhagen, Universitetsparken 5,
DK-2100 Copenhagen {\O}, Denmark. \\ E-mail: michael@math.ku.dk

\end{abstract}


\sect{Introduction}

In this paper we propose a simple and generally applicable method for simulation of
a multi-variate diffusion bridge. The main  
motivation is that simulation of diffusion bridges plays
a fundamental role in simulation-based likelihood inference (including
Bayesian inference) for discretely sampled diffusion processes and other 
diffusion-type processes like stochastic volatility models.

Our approach is based on the following simple construction of a
process that starts from $a$ at time zero and at time $T$ ends in $b$,
where $a$ and $b$ are given points in the state space. One diffusion process,
$X^{(1)}_t$, is started from the point $a$, while another diffusion, $X^{(2)}_t$ is started from the
point $b$. The time of the second diffusion is reversed, so that the
time starts at $T$ and goes downwards to zero, and the dynamics of
$X^{(2)}_t$ is chosen such that the time reversed diffusion
$X^{(2)}_{T-t}$ has the same dynamics as $X^{(1)}_t$. Suppose there is a time
point $\tau \in [0,T]$ at which $X^{(1)}_\tau = X^{(2)}_{T-\tau}$. Then the process that is equal to
$X^{(1)}_t$ for $t \in [0,\tau]$ and for $t \in [\tau,T]$ equals $X^{(2)}_{T-t}$ is obviously a
process that starts at $a$ and ends at $b$. 
If the two diffusion processes $X^{(1)}$ and $X^{(2)}$ are
independent, then the probability that $X^{(1)}_t$ and $X^{(2)}_{T-t}$ meet
at the same time point in $[0,T]$ is zero for dimensions larger than
one. However, if they are suitably dependent, then the probability can
be made positive and will often go to one as $T$ tends to infinity. This can
be obtained by applying classical coupling methods. 

The new method is a generalization of the one-dimensional diffusion
bridge simulation method proposed by \cite{bladtsorensen}, where the
two diffusions $X^{(1)}$ and $X^{(2)}$ were independent. The generalization
is far from straightforward. For ergodic
one-dimensional diffusions, the probability that two independent
diffusions intersect goes to one as $T \rightarrow \infty$. The
application of coupling methods is a breakthrough that allows the
generalization to multivariate diffusions. Moreover, the coupling methods
also improve the simulation of one-dimensional diffusions because they 
increase the probability of intersection in $[0,T]$ and thus
improve the computational efficiency of the method. 

\begin{figure}
\begin{center}
\includegraphics[width=6cm]{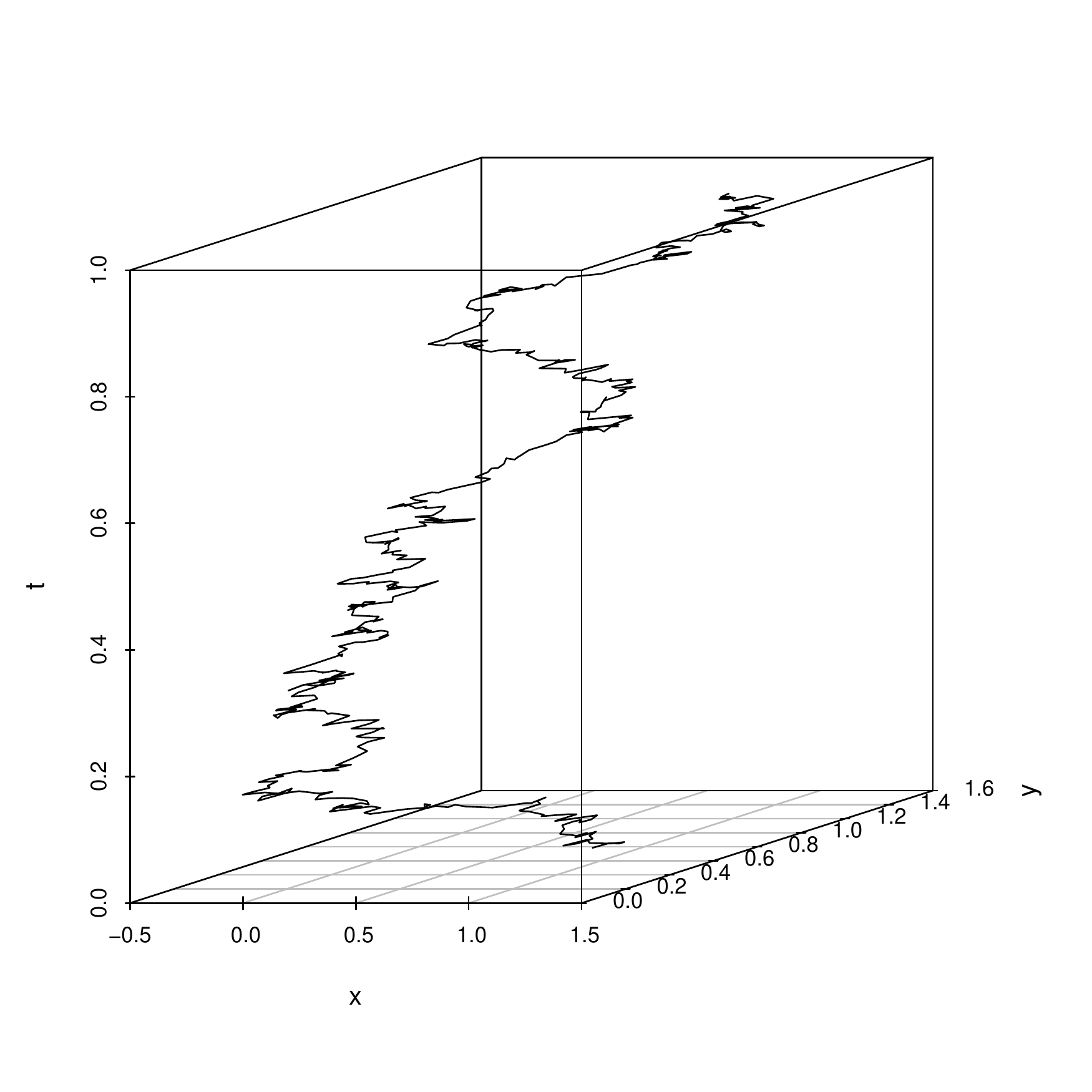} 
\caption{\label{3D}
Simulated sample path of an Ornstein-Uhlenbeck bridge from
(0.785,0.785) to (1.091,1.091). Time is in the vertical direction.}
\end{center}
\end{figure}

Conditional on the event that the two processes meet at a time point
in $[0,T]$, we show that the process constructed as described above is
an approximation to a diffusion bridge 
between the two points. A simple rejection sampler is obtained by
repeatedly simulating the two dependent diffusions until they hit each 
other. The 
diffusions can be simulated by means of simple procedures like the
Euler or the Milstein scheme, see \cite{KlPl}, so the new method is easy to
implement for likelihood inference for discretely sampled diffusion
processes. The approximate diffusion bridge produced by the rejection
sampler can be used as proposal for MCMC-algorithms that have an 
exact diffusion bridge as target distribution. We present a pseudo-marginal 
Metropolis-Hastings algorithm (in the sense of \cite{andrieuroberts}) and 
a new MCMC algorithm that is easier to implement, but typically has a 
larger rejection probability. An example of a diffusion bridge simulated by 
our new method is shown in Figure \ref{3D}.

Diffusion bridge simulation is a highly non-trivial problem that has
been investigated actively over the last 10 - 15 years. A lucid exposition
of the problems and the state-of-the-art can be found in \cite{OPGR}.
Before the paper by \cite{bladtsorensen}, it was thought impossible to
simulate diffusion bridges by means of simple procedures, because a
rejection sampler that tries to hit the prescribed end-point for the
bridge will have a prohibitively large 
rejection probability. The rejection sampler presented in this paper 
has an acceptable rejection probability because what must be hit 
is a sample path rather than a point and because the coupling methods
make the two diffusions tend to meet. The first diffusion bridge simulation 
methods in the literature were based on the Metropolis-Hastings algorithm
with a proposal distribution given by a diffusion process that was
forced by its drift to go
from $a$ to $b$, see e.g.\ \cite{roberts} or \cite{durhamgallant}. 
Later \cite{Beskos:2007,beskos7}
developed algorithms for exact simulation of diffusion bridges. These
are cleverly designed rejection sampling algorithms that use simulations of
Brownian bridges, which can easily be simulated. Under strong
boundedness conditions the algorithm is relatively simple, whereas it is more
complex under weaker condition. \cite{linchenmykland} proposed a
sequential Monte Carlo method for simulating diffusion bridges with a
resampling scheme guided by the empirical distribution of backward
paths. The spirit of this approach has similarities to the methods
in \cite{bladtsorensen} and in this paper.

An advantage of our new method is that the same simple algorithm can
be used for all ergodic diffusions, and that it is easy to understand and to
implement. More importantly, the method does not require
that the diffusion can be transformed into one with diffusion matrix equal
to the identity matrix. Such a transformation, often referred to as the Lamperti
transformation, exists for only a small subclass
of the multi-variate diffusions, and even when it exists, the transformation is
rarely in closed form. A Lamperti transformation is required for the exact algorithms of
\cite{Beskos:2007,beskos7}. Another important advantage
is that our method works particularly well for long time
intervals. The computational complexity is linear in the length of the
time interval where the diffusion bridge is defined. This was
illustrated in a simulation study in \cite{bladtsorensen}, where the
computer time increased linearly with the interval length for our
method, while it grew at least exponentially with the interval length for
the exact EA algorithms of \cite{Beskos:2007}. The latter finding is
not surprising, because in the fundamental EA1 algorithm the
acceptance probability is of the order $e^{-cT}$, where $T$ is the
length of the interval. Thus the EA algorithm
is in practice likely not to work for long time intervals. It follows from results
in this paper that under conditions given in the literature on
coupling of diffusion processes (see \cite{ChenLi}), the approximate method proposed here
simulates an essentially exact diffusion bridge in long time intervals (apart from the discretization
error). This literature also gives conditions ensuring that the distribution of the simulated process
goes to that of a diffusion bridge exponentially fast as a function of
the interval length. Thus the proposed method usefully supplements
previously published methods both because it works particularly well
for long time intervals, where the other methods tend not to work, and
because it works for diffusions without a Lamperti transformation. It is
worth noting that simulation-based likelihood inference for discretely
sampled diffusions is mainly important for long time intervals,
because for short time intervals several simpler methods provide
highly efficient estimators, see the following discussion.

The main challenge in likelihood based inference for diffusion models is that
the transition density, and hence the likelihood function, is not
explicitly available and must therefore be approximated. When the
sampling frequency is relatively high, rather crude approximations to the likelihood
functions, like those in \cite{ozaki}, \cite{bollerslev},
\cite{bmb&ms1} and \cite{kessler97},  
give estimators with a high efficiency, see \cite{ms10}. When the
interval between the observation times is relatively long, more
accurate approximations to the transition density are needed. One approach
is numerical approximations, either by solving the Kolmogorov  PDE
numerically, e.g.\ \cite{poulsen} and \cite{wood&trees}, or by
expansions, e.g.\ 
\cite{ait-sahalia,ait-sahalia2} and
\cite{julie}. Alternatively, likelihood inference can be based on
simulations, an approach that goes back to the seminal paper by
\cite{pedersen}. The inference problem can be
viewed as an incomplete data problem. If the diffusion 
process had been observed continuously, the likelihood functions would
be explicitly given by the Girsanov formula. However, the process has
been observed at discrete time points only, and the continuous-time
paths between the observation points can be considered as missing
data. This way of viewing the problem, which goes back to
\cite{castellezmirou}, makes it natural to apply either the
EM-algorithm or the Gibbs sampler. To do so the missing continuous
paths between the observations must be simulated conditional on the
observations, which by the Markov property is exactly simulation of
diffusion bridges. It was a significant break-through when this was
simultaneously realized by several authors, see \cite{roberts},
\cite{elerian}, \cite{eraker}, and \cite{durhamgallant}, and  
approaches based on bridge simulation has since been used by several 
authors including 
\cite{golightly&wilkinson1,golightly&wilkinson3,golightly&wilkinson2},
\cite{Beskos:2006}, \cite{delyon&hu}, \cite{beskos9}, and
\cite{linchenmykland}.

Diffusion bridge simulation is also crucial to simulation-based
inference for other types of diffusion process data than discrete time
observations. \cite{chibpittshephard} presented a general approach to
simulation-based Bayesian inference for diffusion models when the data
are discrete time observations of rather general, and possibly random,
functionals of the continuous sample path, see also
\cite{golightly&wilkinson4}. This approach covers
for instance diffusions observed discretely with measurement error and
discretely sampled stochastic volatility models. Also in this case
diffusion bridge simulation is crucial. \cite{fernando} presented an
EM-algorithm for integrated diffusions observed discretely with
measurement error based on the ideas in \cite{chibpittshephard} and
the bridge simulation method proposed in \cite{bladtsorensen}.

The paper is organized as follows. In Section 2 we first review some
necessary results on coupling methods and time-reversal for diffusion
processes and prove some preliminary results. Then we present the new
approximate bridge simulation method and show in what sense it
approximates a diffusion bridge. The approximate bridges are then used
as proposal in two MCMC-algorithms that have an exact diffusion bridge 
as target distribution. Finally, we discuss how the new method
improves simulation of 
one-dimensional diffusions and solve some implementation problems. In 
particular, we give criteria to determine whether two diffusions 
simulated at discrete time points have met between two time points. In 
Section 3 the approximate and exact bridge simulation methods are compared 
to the (known) exact distribution of the multivariate Ornstein--Uhlenbeck 
bridge. The study indicates that our approximate method provides a very 
accurate approximation to the distribution of a diffusion bridge, except 
for bridges that are unlikely to occur in discretely sampled data. 
Even for extremely unlikely bridges, the approximate 
method works surprisingly well for some coupling methods. 
In Section 4 we illustrate the usefulness of our
method to inference for discretely observed diffusions by considering briefly
Bayesian estimation for the multivariate hyperbolic diffusion.
The proofs are collected in Section 5.   

\vspace{10mm}

\sect{Diffusion bridge simulation}

Let $X = \{ X_t \}_{t\geq 0}$ be a $d$-dimensional diffusion with
state space $D$ given by the stochastic differential equation  
\be
\label{basicsde}
dX_t = \alpha(X_t)dt + \sigma (X_t)dW_t ,  
\ee
where $W$ is a $d$-dimensional Wiener process, and where the
coefficients $\alpha$ (a function $D \mapsto \R^d$) and $\sigma$ (a
$d \times d$-matrix of continuous functions defined on $D$) are
sufficiently regular to ensure that the 
equation has a unique strong solution that is a strong Markov process.
We will assume that the diffusion defined by (\ref{basicsde})
is ergodic with invariant probability density function $\nu$ (w.r.t.\
Lebesgue measure on $D$). It is assumed that $\sigma(x)$ is invertible
for all $x \in D$. Define $V(x)=\sigma(x)\sigma(x)^T$.
We denote the transition density of $X$ by $p_t(x,y)$. Specifically, 
the conditional density of $X_{s+t}$ given $X_s=x$ is $y \mapsto
p_t(x,y)$.

Let $a$ and $b$ be given points in $D$. We present  
a method for simulating a sample path of $X$ in $[0,T]$ such that
$X_0=a$ and $X_T=b$.
A solution of (\ref{basicsde}) in the interval $[0,T]$ such that 
$X_{0}=a$ and $X_{T}=b$ will in the following be called an 
$(a,b,T)$-bridge. The approximate bridge construction goes as
follows. First the time-reversed version, $X^*_t$, of (\ref{basicsde}) is
simulated in $[0,T]$ starting at $b$. Then a solution, $X_t$, to
(\ref{basicsde}) is simulated starting at $a$ and dependent
on $X^*$ in such a way that $X_t$ and $X^*_{T-t}$ tend to intersect
at some (random) time point $\tau \in [0,T]$. If the two processes meet
at time $\tau$, then the approximate bridge is the equal to
$X_t$ for $t \in [0,\tau]$, and equal to $X^*_{T-t}$ for $t \in
[\tau,T]$. This approximate bridge is used as a proposal for a MCMC
algorithm that has the exact $(a,b,T)$-bridge as target
distribution. The required dependence between $X^*$ and $X$ is
obtained by applying classical coupling methods for diffusions.

\vspace{3mm}

\subsection{Coupling and time-reversal for multivariate diffusions}

We present our algorithm
for a class of coupling methods that includes the coupling by
reflection method of \cite{lindvallrogers} and the coupling by
projection by \cite{ChenLi}. Other coupling methods 
(see e.g.\ \cite{ChenLi}) can be used similarly, provided that they 
couple before time $T$ with a probability that is not too small.
We begin by briefly presenting the class of coupling methods. Then we derive
a few results that we need in order to construct diffusion bridges. 

Suppose $X$ solves (\ref{basicsde}) with initial value
$X_0=x_0$. Define another diffusion process $X'$ as the solution to
\be
\label{couplesde}
dX'_t = \alpha(X'_t)dt + \sigma (X'_t)dW'_t , \ \ \ X'_0 = x'_0 
\ee
with the Wiener process
\[
dW'_t = \left\{ I-(1-\gamma)\Pi(X_t,X'_t) \right\} {\cal O}(X_t,X'_t)dW_t 
+ \sqrt{1-\gamma^2} u(X_t,X'_t)dU_t.
\]
Here $\gamma \in [-1,1)$, $U$ is a univariate standard Wiener process
independent of $W$, $I$ is the $d$-dimensional identity matrix and
\be
\label{V}
\Pi(x,x')=u(x,x')u(x,x')^T,
\ee
where $T$ denotes transposition, and $u(x,x')$ is the unit vector such 
that $\sigma(x')u(x,x')$ points in the direction $x-x'$, i.e.
\[
u(x,x') = \frac{\sigma(x')^{-1}(x-x')}{|\sigma(x')^{-1}(x-x')|}.
\]
Finally, ${\cal O}(X_t,X'_t)$ is an orthogonal matrix that in some
cases is needed to ensures that the law of $(X,X')$ does not depend on
the particular choice of $\sigma$, but only on the law of the solution
$X$ (the law of $X$ depends only on $\sigma \sigma^T$, so the same law
can be obtained for many different choices of $\sigma$). For many stochastic 
differential equations, indeed for all those we consider in this paper, we
have ${\cal O} = I$, but there may be cases  where the choice of
$\sigma$ introduces a rotation which should be counterbalanced.  
In fact, as long as $\sigma$ satisfies the smoothness requirements of 
\cite{lindvallrogers} and \cite{ChenLi} our method can still work if we choose ${\cal O} = I$
for reasons of computational efficiency. However, it may take many
more attempts to achieve the successful coupling from which we construct our bridge.
The matrix ${\cal O}(X_t,X'_t)$ should be chosen to be the closest
orthogonal matrix to $\sigma (X_t)^{T}\sigma(X'_t)$ 
in the Frobenious norm. That is ${\cal O} (X_t, X_t') = AB^T$ where 
$\sigma (X_t)^{T}\sigma(X'_t) = A\Sigma B^T$ is the singular value decomposition.

The matrix $\Pi(x,x')$ is the orthogonal projection onto the one-dimensional 
subspace generated by the vector $u(x,x')$, while $I-\Pi(x,x')$ is 
projection onto the plane orthogonal to the vector $u(x,x')$. Using
this geometric interpretation, it is 
not difficult to see that the quadratic variation of $W'$ equals $tI$
implying that $W'$ is a Wiener process. 

Consider the case where ${\cal O} = I$. In the plane orthogonal to the
vector $u(X_t,X'_t)$, the increment of the Wiener process $W'$ is
equal to the increment of $W$. In the direction $u(X_t,X'_t)$, the
increment of $W'$ is equal to minus the increment of $W$ in the same
direction (i.e.\ $\Pi(X_t,X'_t)dW_t$) if $\gamma=-1$ ({\it method of 
reflection}). Otherwise, the increment of $W'$ in the direction 
$u(X_t,X'_t)$ is the sum of $\gamma \Pi(X_t,X'_t)dW_t$ and
$\sqrt{1-\gamma^2} u(X_t,X'_t)dU_t$. In particular if $\gamma = 0$, the 
increment of $W'$ in the direction $u(X_t,X'_t)$ it is equal to the
increment of the independent Wiener process $U$ on the subspace
generated by $u(X_t,X'_t)$ ({\it method of projection}). For $\gamma =
-1$, $dW'_t = H(X_t,X'_t)) dW_t$, where the matrix $$H(x,x')=I-2\Pi(x,x')$$
is reflection in the plane orthogonal to the vector
$u(x,x')$. It is therefore symmetric and orthonormal. We
do not consider the case $\gamma = 1$ where (for ${\cal O} = I$) the
two diffusions are driven by the same Wiener process and will not meet.

The squared diffusion coefficient of the $2d$-dimensional diffusion
$(X_t,X'_t)$ is
\[
\left( \begin{array}{cc} \sigma(x)\sigma(x)^T & \sigma(x){\cal O}(x,x')^T
\left\{ I - (1-\gamma)\Pi(x,x')) \right\} \sigma(x')^T \\ & \\ \sigma(x')
\left\{ I - (1-\gamma)\Pi(x,x')) \right\} {\cal O}(x,x') \sigma(x)^T
  & \sigma(x')\sigma(x')^T \end{array} \right),
\]
which is of the general form treated in \cite{ChenLi}.

Define the stopping time
\be
\label{couplingtime}
\tau = \inf \{ t>0 \, | \, X_t = X'_t \}.
\ee
\cite{lindvallrogers}, \cite{ChenLi} and others have given conditions on the
coefficients $\alpha$ and $\sigma$ ensuring that $P(\tau < \infty)
=1$. This is not really necessary for our bridge simulation method to
work. The method is a rejection sampler with rejection probability
$P(\tau > T)$, so we just need that this probability is not too
large. If $P(\tau < \infty) =1$ that is certainly the case if $T$ is
sufficiently large. \cite{ChenLi}) gave results on the rate of
convergence of $P(\tau > T)$ to zero as $T \rightarrow \infty$, in
particular conditions ensuring geometrically fast convergence. As
these conditions are somewhat technical and unnecessarily restrictive
for our application, they are not stated here.

\begin{lem}\label{lemma:samplepaths}
The sample path of $X'$ in $[0,t]$ is a
function of the sample path of $X$ in $[0,t]$, the initial value
$x'_0$, and the sample path of the one-dimensional Wiener process 
$U$ in $[0,t]$
\[
\{ X'_s \}_{0 \leq s \leq t} = {\cal K}_t \left( x'_0, \{ X_s \}_{0 
\leq s  \leq t}, \{ U_s \}_{0 \leq s  \leq t} \right).  
\]
Specifically,
\bea
\lefteqn{X'_s = x'_0 + \int_0^s \alpha(X'_u)du + \int_0^s \sigma(X'_u) 
\{ I-(1-\gamma)\Pi(X_t,X'_t)\} {\cal O}(X_t,X'_t) \sigma (X_u)^{-1}
dX_u} \nonumber  \\ &&  \mbox{} \hspace{15mm}
- \int_0^s \sigma(X'_u) \{ I-(1-\gamma)\Pi(X_t,X'_t)\} {\cal O}(X_t,X'_t)
\sigma (X_u)^{-1} \alpha (X_u)du \nonumber  \\ &&  \mbox{}  \hspace{15mm}
+ \sqrt{1-\gamma^2} \int_0^s \sigma(X'_u) u(X_t,X'_t) dU_u.
\label{prime-dif}
\eea
Similarly, the sample path of $X$ in $[0,t]$ is a function of the initial 
value $x_0$, the sample paths of $X'$ in $[0,t]$, and the sample path
in $[0,t]$ of a standard univariate Wiener process, $U'$, 
independent of $X'$
\[
\{ X_s \}_{0 \leq s \leq t} = \tilde {\cal K}_t 
\left( x_0, \{ X'_s \}_{0 \leq s  \leq t}, \{ U'_s \}_{0 \leq s  
\leq t} \right). 
\]
Specifically,
\bea
X_s &=& x_0 + \int_0^s \alpha(X_u)du + \int_0^s \sigma(X_u) 
{\cal O}(X_u,X'_u)^T \{ I-(1-\gamma)\Pi(X_u,X'_u)\} \sigma (X'_u)^{-1}
dX'_u \nonumber \\ &&  \mbox{} \hspace{15mm}
- \int_0^s \sigma(X_u) {\cal O}(X_u,X'_u)^T \{ I
-(1-\gamma)\Pi(X_u,X'_u)\} \sigma (X'_u)^{-1} \alpha (X'_u)du
 \nonumber \\ &&  \mbox{}  \hspace{15mm} + \sqrt{1-\gamma^2} 
\int_0^s \sigma(X_u) {\cal O}(X_u,X'_u)^T u(X_u,X'_u) dU'_u.
\label{stjerne-dif}
\eea
\end{lem}

Proofs of this and other results are given in Section \ref{proofs}.
Note that for $\gamma = -1$, $\tilde {\cal K}_t$ does not depend on $U'$.

Before we can formulate the main theorem on our method, we need to
review some well-known results on time-reversal of multivariate
diffusions. We have assumed that the diffusion $X$ defined by
(\ref{basicsde}) is ergodic with invariant probability density
function $\nu$. Hence a stationary version of $X$ exists. If 
the time is reversed for this stationary process, we obtain another
stationary diffusion process $X^*$. By Theorem 2.3 in \cite{nualart},
the time-reversed diffusion $X^*$ solves the stochastic differential
equation
\be
\label{reversesde}
dX^*_t = \alpha^*(X^*_t)dt + \sigma (X^*_t)dW_t
\ee
where
\be
\label{reversedrift}
\alpha_i^*(x) =  -\alpha_i(x) + \nu(x)^{-1}\sum_{j=1}^d \partial_{x_j} 
\left( \nu(x) V(x)_{ij} \right), \ \ i = 1, \ldots, d, 
\ee
provided that 
\be
\label{sdebetingelse}
\int_D \left| \sum_{j=1}^d \partial_{x_j} \left( \nu(x) V_{ij}(x) \right) 
\right| dx < \infty, \ \ i = 1, \ldots , d.
\ee
We assume that the local integrability condition (\ref{sdebetingelse})
is satisfied. Conditions ensuring this are discussed in
\cite{nualart}, where also a similar result under the
local Lipschitz condition is given. The condition (\ref{sdebetingelse})
is satisfied if the two coefficients are twice continuously differentiable 
on $D$, and if there exists $\epsilon > 0$ such that $V \geq
\epsilon I$. Alternative conditions
can be found in \cite{pardoux}. In cases where the transition density
is not differentiable with respect  to $x$, the partial derivative in
the formula for the drift are in the distributional sense.

Let $p^*_t(x,y)$ denote the
transition density of the solution to (\ref{reversesde}). If $X$ is
the stationary version of the solution to (\ref{basicsde}) and $X^*$
is the time-reversed stationary diffusion, then the distribution of
$(X_s,X_{s+t})$ equals the distribution of $(X^*_{s+t},X^*_s)$. Hence
\be
\label{balance}
p_t(x,y)\nu(x) = p^*_t(y,x)\nu(y),
\ee

\vspace{3mm}

\subsection{Approximate bridge simulation}
\label{approxbridge}

In this subsection we present the mathematical results on which our
algorithm to approximately simulate a diffusion bridge is based and
describe how the results can be used to construct the algorithm. 
Detailed implementation questions are discussed in a later section. 

\begin{theo}
\label{theorem1}
Suppose $X$ solves (\ref{basicsde}) for $t \in [0,T]$ with the initial
condition $X_0 \sim \nu$, where $\nu$ is the invariant probability
measure. Let $X'$ be the corresponding solution to
(\ref{couplesde}) with initial condition $X'_0=a$, i.e.\ 
$\{ X'_t \}_{0 \leq t \leq T} = {\cal K}_T \left( a, \{ X_s \}_{0 
\leq t  \leq T}, \{ U_t \}_{0 \leq t  \leq T} \right)$, where $U$ is a
standard Wiener process independent of $X$. Define a process by
\[ 
Z_t=\left\{ 
\begin{array}{ll}
X'_t & \mbox{\rm if } \ 0\leq t \leq \tau \\ & \\
X_{t} & \mbox{\rm if } \ \tau < t \leq T, 
\end{array}
\right.
\]
where $\tau$ is given by (\ref{couplingtime}).

Then the distribution of $\{ Z_t \}_{0\leq t \leq T}$ conditional on the
events $\tau \leq T$ and $X_T=b$  equals the distribution of a 
$(a,b,T)$-bridge, $B$, conditional on the event that the bridge is
hit by the process $\tilde {\cal K}_T \left( A , \{ B_t \}_{0 
\leq t  \leq T}, \{ U'_t \}_{0 \leq t  \leq T} \right)$.
Here $A$ is a $d$-dimensional random variable with density function
$p^*_T(b,\cdot)$ given by (\ref{balance}), $U'$ is a standard univariate
Wiener process, and $A$, $U'$ and $B$ are independent.
\end{theo}

\noindent
We refer to the process  $\tilde {\cal K}_T \left( A , 
\{ B_t \}_{0 \leq t  \leq T}, \{ U'_t \}_{0 \leq t  \leq T} \right)$
as the $p^*_T(b)$-{\it diffusion associated with} $B$. This process
plays an important role not only in
the characterization of the distribution of the approximate diffusion
bridge $Z$, but also in the method for simulating exact diffusion
bridges presented in the next subsection. When $\gamma = -1$, 
${\cal K}_T$ does not depend on $U$
and $\tilde {\cal K}_T$ does not depend on $U'$.

A sample path of $X$ with $X_0 \sim \nu$ conditional on $X_T = b$
can easily be obtained by using the result of the following 
lemma on the distribution of a time-reversed diffusion started at the
point $b$. 

\begin{lem}\label{lemma:fundamental}
Suppose $X$ is ergodic with invariant probability density
function $\nu$, and
let $X^*$ be a solution to (\ref{reversesde}) with initial condition
$X^*_0 = b$. Define the time-reversed process $\bar X_t = X^*_{T-t}$,
$0 \leq t \leq T$. The process $\{ \bar X_t  \}$ and the conditional
process $\{ X_t \}$ given that $X_T=b$ have the same transition densities
\be
\label{bridgetransition}
q(x,s,y,t)= \frac{p_{t-s}(x,y)p_{T-t}(y,b)}{p_{T-s}(x,b)}
= \frac{p^*_{t-s}(y,x)p^*_{T-t}(b,y)}{p^*_{T-s}(b,x)}, 
\hspace{5mm} s < t < T.
\ee
The distribution of $\{
\bar X_t  \}$ is equal to the distribution of the process
$\{ X_t \}$ with $X_0 \sim \nu$ conditional on $X_T=b$.
\end{lem}

Based on Theorem \ref{theorem1} and Lemma \ref{lemma:fundamental} we
can now propose an algorithm to simulate an approximate diffusion
bridge in the interval $[0,T]$. Use any of the several methods
available (see e.g.\ \cite{KlPl}) to simulate the diffusion $X^*$
given by (\ref{reversesde}) with $X^*_0=b$. If the diffusion
given by (\ref{basicsde}) is time-reversible, then the stochastic
differential equation for $X^*$ is simply (\ref{basicsde}). To
simplify the exposition, we assume that $X^*$ has been simulated by
means of the Euler-scheme with step size $\delta = T/N$. Let 
$Y^*_{\delta i}$, $i=0,1, \ldots, N$, denote the simulated values of
the process, while $\Delta W_i = W_{\delta i} - W_{\delta (i-1)}$, 
$i=1, \ldots, N$, denote the simulated increments of the driving
$d$-dimensional Wiener process, i.e.\ $Y^*_0 = b$ and
\[
Y^*_{\delta i} = Y^*_{\delta (i-1)} + \alpha^*(Y^*_{\delta (i-1)})\delta 
+ \sigma (Y^*_{\delta (i-1)}) \Delta W_i,
\]
$i=1,\ldots,N$. The increments of the Wiener process that drives the
time-reversed version of $Y^*$ (i.e.\ $Y^{*\mbox{rev}}_{\delta i} =
Y^*_{\delta (N-i)}$) are
\be
\label{Wrev}
\Delta W^{\mbox{rev}}_i = \sigma (Y^*_{\delta (N-i+1)})^{-1}\left(
Y^*_{\delta (N-i)} - Y^*_{\delta (N-i+1)} - \alpha (Y^*_{\delta (N-i+1)})
\delta \right).
\ee
For fine discretizations $\Delta W^{\mbox{rev}}_i \approx - \Delta
W_{N-i+1}$.

The discretized sample path of $X'$ is a function of the simulated
process $Y^*$ and (except in the case of the method of reflection, 
$\gamma = -1$) an independent one-dimensional standard Wiener
process $B$, the increments of which we denote by $\Delta B_i =
B_{\delta i} - B_{\delta (i-1)}$, $i=1, \ldots, N$. If we
denote the simulated values of $X'$ by $Y'_{\delta i}$, $i=0,1,
\ldots, N$, we have that $Y'_0 = a$ and
\be
\label{simY'}
Y'_{\delta i} = Y'_{\delta (i-1)} + \alpha (Y'_{\delta (i-1)})\delta 
+\sigma (Y'_{\delta (i-1)}) \Delta W'_i,
\ee
$i=1,\ldots,N$, where
\bea
\label{simW'}
\Delta W'_i &=& \left\{I - (1-\gamma)\Pi(Y^*_{\delta (N-i+1)},
Y'_{\delta (i-1)})\right\} {\cal O}(Y^*_{\delta (N-i+1)},Y'_{\delta (i-1)})
\Delta W^{\mbox{rev}}_i \\ && \mbox{} \hspace{20mm} 
+ \sqrt{1-\gamma^2}u(Y^*_{\delta (N-i+1)},Y'_{\delta (i-1)}) \Delta B_i.
\nonumber
\eea

A simulation of an approximation in the sense of Theorem
\ref{theorem1} to a $(a,b,T)$-bridge is obtained by 
rejection sampling. Keep simulating independent copies of $Y^*$
and $B$
until there is an $i$ such that $Y^*_{\delta (T-i)}$ and $Y'_{\delta i}$ 
are sufficiently close that we can safely assume that coupling happens
in the time interval $[\delta i, \delta (i+1)]$. We discuss
the problem of deciding whether coupling has happened or not
in more detail in Subsection \ref{imp}. Once coupling has been
obtained (in the interval $[\delta i, \delta (i+1)]$), put $\rho =
i+1$ and define 
\be
\label{bridgesim}
Z_{\delta i}=\left\{ 
\begin{array}{ll}
Y'_{\delta i} & \mbox{\rm for } \ i = 0,1, \ldots, \rho-1 \\ & \\
Y^*_{\delta (N-i)} & \mbox{\rm for } \ i=\rho, \ldots N, 
\end{array}
\right.
\ee
On top of the usual influence of the step size
$\delta$ on the quality of the individual simulated trajectories, the
step size also controls the probability that a trajectory crossing is
not detected. Therefore, it is advisable to choose $\delta$ smaller
than in usual simulation of diffusion sample paths. 
Another problem is that the method will only work, if $P(\tau \leq T)$
is not too small. This problem was considered in
\cite{lindvallrogers} and \cite{ChenLi}.

The results of Lemma \ref{lemma:fundamental} and Theorem \ref{theorem1}, 
and hence the algorithm, simplify if the diffusion process is time
reversible in the sense that $p^*_t(x,y)=p_t(x,y)$, or equivalently 
$p_t(x,y)\nu(x) = p_t(y,x)\nu(y)$. Diffusions with the latter property
are called $\nu$-symmetric, see  \cite{kent}. By equating $\alpha$ to
the reverse drift $\alpha^*$ given by (\ref{reversedrift}), 
it follows that the diffusion given by (\ref{basicsde}) is time-reversible if
\be
\label{timereversible}
\alpha_i(x) \nu(x) =  \mbox{\small $\frac12$} \sum_{j=1}^p 
\partial_{x_j} \left( \nu(x) V(x)_{ij} \right), \ \ i = 1, \ldots, d. 
\ee
When $V(x)$ is a diagonal matrix, these equations simplify to
\be
\label{timereversible_simple}
\alpha_i(x) \nu(x) =  \mbox{\small $\frac12$} \partial_{x_i} 
\left( \nu(x) V(x)_{ii} \right), \ \ i = 1, \ldots, d. 
\ee
  
\vspace{3mm}
  
\subsection{Exact bridge simulation}
\label{exactsimulation}

The algorithm presented in the previous section produces only
approximate diffusion bridges. The simulations in Section \ref{simul}
indicate that the approximation is usually good, and for certain
coupling methods can be even very good. In order to produce
exact diffusion bridges, this section presents two MCMC methods that use
the approximate bridges as proposals and have an exact diffusion bridge
as target distribution. 

First we investigate how the distribution of the approximate diffusion
bridge is related to the distribution of the exact diffusion bridge.
Diffusions, diffusion bridges and the approximate diffusion bridge $Z$
are elements of the canonical space, $C_T$, of $\R^d$-valued
continuous functions defined on the time interval $[0,T]$. Each of
these processes induce a probability measure on the usual
$\sigma$-algebra generated by the cylinder sets. Let $f_{br}$ denote the
Radon-Nikodym derivative of the distribution of the
$(a,b,T)$-diffusion bridge with respect to a dominating measure. The
diffusion bridge solves a stochastic differential equation with the
same diffusion coefficient as in (\ref{basicsde}), see e.g.\ (4.4) in
\cite{OPGR}, so the density $f_{br}$ is given by Girsanov's
theorem. Since the drift for the bridge is unbounded at the end point,
one has to choose the dominating measure carefully: it must correspond
to another bridge, see \cite{OPGR}, p.\ 322, and \cite{delyon&hu}. 
Similarly let $f_a$ denote the density of the distribution of the
approximate bridge $Z$. 

We know from Theorem \ref{theorem1} that the relation between the
distribution of a diffusion bridge $B$
and the approximation $Z$ involves the $p^*_T(b)$-diffusion associated
with $B$, i.e.\ the process $\tilde {\cal K}_T \left( A, B, U \right)$. 
Here $A$ is a $d$-dimensional random variable with density function 
$p^*_T(b,\cdot)$, and $U$ is a one-dimensional standard Wiener process,
where $A$, $B$ and $U$ are independent.

For any $x \in C_T$, let $M_x$ be the set of functions $y \in
C_T$ that intersect $x$. Specifically,
\[
M_x = \{ y \in C_T \, | \, \mbox{gr}(y) \cap \mbox{gr}(x) 
\neq \emptyset  \},
\]
where $\mbox{gr}(x) = \{ (t,x_t) \, | \, t \in [0,T]  \}$. With these
definitions, the relation between the distribution of the approximate
bridge $Z$ and the exact $(a,b,T)$-diffusion bridge is given by the
following corollary.

\begin{cor}
The density of the approximate bridge $Z$ is given by
\label{cor3}
\be
\label{z-density}
f_a(z) = f_{br}(z) \pi_T (z)/\pi_T,
\ee
where $f_{br}$ is the density of an exact diffusion bridge, and
\be
\label{pi1}
\pi_T(x) = P(\tilde {\cal K}_T \left( A, x, U \right) \in M_x),
\hspace{5mm}
\pi_T = P((B,\tilde {\cal K}_T \left( A, B, U \right)) \in M),
\ee
where $M = \{ (x,y) \in C_T \times C_T \, | \, y \in M_x \}$, 
$B$, $A$ and $U$ are independent, $B$ is a $(a,b,T)$-diffusion
bridge,  and $A$ is a $d$-dimensional random variable with density
function $p^*_T(b,\cdot)$, and $U$ is a standard univariate Wiener process.
\end{cor}

Clearly, $\pi_T(x)$ is the
probability that a trajectory $x$ is hit by the $p^*_T(b)$-diffusion
associated with $x$, while $\pi_T$ is probability that an $(a,b,T)$-bridge
is hit by its associated $p^*_T(b)$-diffusion. 
Equation (\ref{z-density}) gives an explicit
expression of the quality of our approximate simulation method, and
more importantly, it can be used to construct MCMC-algorithms that
have the exact distribution of an  $(a,b,T)$-diffusion bridge as
target distribution.

Simulation of $p^*_T(b)$-diffusions associated to a given simulated sample 
path $Z$ of an approximate $(a,b,T)$-diffusion bridge is crucial to the 
following MCMC-algorithms, so we explain in detail how a sample path
of this process can be simulated. We denote the simulated values by
$\tilde Y_{\delta i}$. First the initial value 
$\tilde Y_0 = A$ with density function $p^*_T(b,\cdot)$ must be
generated. Usually the transition density $p^*$ of the time-reversed
diffusion is not explicitly known, but a value of $\tilde Y_0$ can easily be
generated by simulating a sample path $X^*$ of the time-reversed diffusion 
given as the 
solution to (\ref{reversesde}) in $[0,T]$ with $X^*_0=b$ (independently of 
$Z$). Then $\tilde Y_0 := X^*_T$ has the density $p^*_T(b,\cdot)$. Using 
the Euler scheme, we obtain a discretized sample path $\tilde Y$ from 
(\ref{stjerne-dif}) as follows:
\bean
\lefteqn{\tilde Y_{\delta i} = \tilde Y_{\delta (i-1)} +
\alpha(\tilde Y_{\delta (i-1)})\delta} \\ && \mbox{} \hspace{8mm}
+ \sigma(\tilde Y_{\delta (i-1)}) {\cal O}(\tilde Y_{\delta (i-1)},
Z_{\delta(i-1)})^T \left\{I - (1-\gamma)\Pi(\tilde Y_{\delta (i-1)},
Z_{\delta(i-1)})\right\} \Delta \tilde W_i \\ && \mbox{} \hspace{8mm} +
\sqrt{1-\gamma^2}u(\tilde Y_{\delta (i-1)},
Z_{\delta(i-1)}) \Delta U_i, \hspace{10mm} i=1, \ldots, T,
\eean
where
\[
\Delta \tilde W_i = \sigma(Z_{\delta (i-1)})^{-1}\left\{ Z_{\delta i}
- Z_{\delta (i-1)} - \alpha (Z_{\delta (i-1)})\delta \right\}
\]
and $\Delta U_i$, $i=1, \ldots, T$, are independent 
$N(0,\delta)$-distributed random
variables (independent of $Z$ and $\tilde Y_0$). The increments
$\Delta \tilde W_i$ were calculated during the simulation of the
approximate bridge $Z$. Specifically, for $i \leq \rho -1$, $\Delta 
\tilde W_i$ equals the Wiener process increment in (\ref{simY'}) given by
(\ref{simW'}). Here $\rho$ is the time where the two process defining
$Z$ meet; cf.\ (\ref{bridgesim}). For $i \geq \rho$, $\Delta \tilde W_i$
equals $\Delta W^{\mbox{rev}}_i$ given by (\ref{Wrev}).

First we present a Metropolis-Hastings algorithm of
the pseudo-marginal type studied in \cite{andrieuroberts} for 
which the proposal is the approximate simulation method with density $f_a$,
and the target distribution is the distribution of an exact diffusion
bridge with density $f_{br}$. A simple MH-algorithm would use a sample 
path $Z^{(i)}$ of an approximate diffusion bridge simulated by one of 
the methods in Subsection \ref{approxbridge} as proposal in the $i$th 
step. The proposed sample path is accepted with probability 
$\alpha(X^{(i-1)}, Z^{(i)}) = \min (1, r(X^{(i-1)}, Z^{(i)}) )$, 
where
\[
r(X^{(i-1)}, Z^{(i)}) = 
\frac{f_{br}(Z^{(i)})f_a(X^{(i-1)})}{f_{br}(X^{(i-1)})f_a(Z^{(i)})} =
\frac{\pi_T(X^{(i-1)})}{\pi_T(Z^{(i)})}.
\]
Here $X^{(i-1)}$ is the sample path from the previous step, and $\pi_T(x)$ 
is the probability given by (\ref{pi1}) that the sample path $x$ is hit 
by the $p^*_T(b)$-diffusion associated with $x$. This MH-algorithm
produces draws of exact diffusion bridges, but the probability $\pi_T(x)$ is
not explicitly known.

As in \cite{bladtsorensen}, exact diffusion bridges can be simulated by
means of a MCMC algorithm of the pseudo-marginal type. 
The basic idea of the pseudo-marginal approach
is to replace the factor in the acceptance ratio, which we cannot calculate,
$f_{br}(x)/f_a(x)$ $= 1/\pi_T(x)$ by an unbiased MCMC estimate. The
beauty of the method is that by including the MCMC draws needed to
estimate $1/\pi_T(x)$ in the MH-Markov chain, the marginal equilibrium
distribution of the bridge draws is exactly $f_{br}$, irrespective of
the randomness of the estimate of $1/\pi_T(x)$. 

For a given sample path $x$ of an approximate diffusion bridge, define 
a random variable $T(x)$ in the following way. Simulate a sequence of 
independent $p^*_T(b)$-diffusion associated with $x$, $\tilde X^{(1)},
\tilde X^{(2)}, \ldots$ until $x$ is intersected by $\tilde
X^{(i)}$, and let $T(x)$ be the index of the first $p^*_T(b)$-diffusion
that hits $x$:
\[
T(x) = \min \{i: \tilde X^{(i)} \in M_x \}.
\]
By results for the geometric distribution $E(T(x)) = 1/\pi_T(x)$,
so if ${\bf T}(x) = (T_1(x),\ldots,T_N(x))$ is a vector of $N$ independent
draws of $T(x)$, then an unbiased and consistent estimator of
$1/\pi_T(x)$ is 
\[
\hat \rho({\bf T}(x)) = \frac1N \sum_{j=1}^N T_j(x).
\]
The pseudo-marginal MH-algorithm goes as follows.
\begin{description}
\item{1.} Simulate an initial approximate diffusion bridge, $X^{(0)}$, by 
by one of the methods in Subsection \ref{approxbridge} and $N$ 
independent (conditionally on $X^{(0)}$) $T(x)$-values,  
${\bf T}^{(0)}(x) = (T_1^{(0)}(x),\ldots,T_N^{(0)}(x))$ with $x=X^{(0)}$, 
and set $i = 1$.
\item{2.} Propose a new sample paths by simulating an approximate
diffusion bridge, $Z^{(i)}$, independently of previous draws (by the 
same method), and simulate $N$
independent (conditionally on $Z^{(i)}$) $T(x)$-values,  
${\bf T}^{(i)} = (T_1^{(i)}(x),\ldots,T_N^{(i)}(x))$ with $x=Z^{(i)}$ 
\item{3.} With probability 
$\min (1, \hat r(X^{(i-1)},{\bf T}^{(i-1)},Z^{(i)},{\bf T}^{(i)}))$, where
\[
\hat r(X^{(i-1)},{\bf T}^{(i-1)}, Z^{(i)},{\bf T}^{(i)}) 
= \frac{\hat \rho({\bf T}^{(i)})}
{\hat \rho({\bf T}^{(i-1)})},
\]
the proposed pair $(Z^{(i)},{\bf T}^{(i)})$ is accepted and $X^{(i)} 
:= Z^{(i)}$. Otherwise $X^{(i)}:= X^{(i-1)}$ and 
${\bf T}^{(i)} := {\bf T}^{(i-1)}$ 
\item{4.} $i := i+1$ and GO TO 2.
\end{description}

By results in \cite{andrieuroberts}, the target distribution of $X$ is
that of an exact diffusion bridge. In fact, since 
\[
\hat r (x^{(1)},{\bf t}^{(1)},x^{(2)},{\bf t}^{(2)}) = 
\frac{f_a(x^{(2)}) {\bf f}_g({\bf t}^{(2)} \, | \, x^{(2)}) 
\hat \rho({\bf t}^{(2)}) 
f_a(x^{(1)}) {\bf f}_g({\bf t}^{(1)} \, | \, x^{(1)})} 
{f_a(x^{(1)}) {\bf f}_g({\bf t}^{(1)} \, | \, x^{(1)})
\hat \rho({\bf t}^{(1)})
f_a(x^{(2)}) {\bf f}_g({\bf t}^{(2)} \, | \, x^{(2)})},
\]
where ${\bf f}_g({\bf t}  \, | \, x)$ is the conditional density of
${\bf T}$ given $X=x$, we see that the density of the target distribution is
\[
p(x,{\bf t}) =  f_a(x) {\bf f}_g({\bf t}  \, | \, x) \hat \rho
({\bf t}) \pi_T
= f_{br}(x) {\bf f}_g({\bf t}  \, | \, x) \hat \rho
({\bf t})\pi_T(x),
\]
where we have used (\ref{z-density}). Since, conditionally on $x$,  $\hat \rho 
({\bf t})$ is an unbiased estimator of $1/\pi_T(x)$,
we find by marginalizing that the target
distribution density of $X$ is $f_{br}$, the density function of an exact
diffusion bridge. 

We also present a simple alternative MCMC algorithm with target 
distribution equal to the exact distribution of a $(a,b,T)$-diffusion 
bridge. The MCMC-algorithm works as follows.
\begin{description}
\item{1.} Simulate an initial approximate diffusion bridge $X^{(0)}$ by 
by one of the methods in Subsection \ref{approxbridge}, and $i := 1$. 
\item{2.} Simulate a $p^*_T(b)$-diffusion $\tilde X^{(i)}$ associated 
with $X^{(i-1)}$.
\item{3.} If $\tilde X^{(i)}$ does not intersect $X^{(i-1)}$, 
then $X^{(i)} := X^{(i-1)}$. Otherwise, simulate a new (independent)
approximate diffusion bridge $X^{(i)}$ (by the same method as in step 1).
\item{4.} $i := i+1$ and GO TO 2.
\end{description}

It is straightforward to check 
that the Markov chain defined in this way satisfies the detailed balance
equation with $f_{br}$ as the stationary density. If $\pi_T(X) > 
\varepsilon > 0$ 
for all potential approximate diffusion bridges $X$, then the Markov chain is
exponentially mixing. To see this, note that as soon as $X^{(i)}$ equals a
new independent approximate diffusion bridge, then $X^{(i)}$ is independent 
of $X^{(0)}$. If $\pi_T(X^{(0)}) > \varepsilon$ with probability one, then 
$P(X^{(i)}=X^{(0)}) < (1-\varepsilon)^i$.

The two MCMC algorithms are probably appropriate for different
applications. For producing a sequence of nearly independent bridges
the pseudo-marginal approach is probably the right choice, while the
alternative MCMC algorithm might be better for calculating expectations.  

To produce diffusion bridges by the proposed MH-algorithm and the
alternative MCMC-algorithm, a number
of sample paths of ordinary diffusions must be simulated. If these
sample paths are simulated by an approximate method, like the Euler
scheme, a small discretization error is introduced. This problem can,
however, in some cases be avoided by using the methods for exactly 
simulating diffusions developed by \cite{Beskos:2007} and 
\cite{beskos7}. This method can be used when a 
multi-dimensional version of the Lamperti transform exists, so that
by this transformation a diffusion can be obtained for which the
diffusion matrix equals the $d$-dimensional identity matrix.
By combining our exact bridge simulation algorithm with exact
diffusion simulation methods, exact diffusion bridges can be
efficiently simulated even in long time intervals.

The computational complexity of both the exact and the approximate 
algorithm is linear in the interval length $T$. The reason is that
the coupling probabilities are non-decreasing functions of the interval
length. Therefore as the interval length increases, the expected
number of rejections when simulating the proposal (the approximate
bridge) is bounded, and the mixing properties of the MCMC-procedure
cannot deteriorate.

\vspace{3mm}

\subsection{One-dimensional diffusions}
\label{1D}

The bridge simulation methods presented in the present paper work in the 
one-dimensional case too and thus generalize the bridge simulation method
proposed in \cite{bladtsorensen}.

For $d=1$ the standard Wiener process driving the process $X'$ given by
(\ref{couplesde}) is simply given by
\[
W'_t = \gamma W_t + \sqrt{1-\gamma^2} U_t, \hspace{5mm} \gamma \in [-1,1),
\]
where $U$ is a standard Wiener process independent of $W$. The bridge 
simulation method in \cite{bladtsorensen} is obtained for $\gamma=0$, 
in which case the diffusion process $X'$ is independent of $X$.

When simulating the $p^*_T(b)$-diffusion associated with $X'$, we need 
to express $X$ in terms of the sample path of $X'$. This is straightforward 
in the one-dimensional case. Define
\[
U'_t = \sqrt{1-\gamma^2}W_t - \gamma U_t.
\]
Then it is easily seen that $U_t$ is a standard Wiener process independent 
of $W'_t$, and that
\[
\gamma W'_t + \sqrt{1-\gamma^2}U'_t = W_t.
\]
Hence
\[
dX_t = \left\{\alpha(X_t)-\gamma \sigma(X_t) \sigma(X'_t)^{-1}
\alpha(X'_t)\right\}dt + \sigma(X_t)\left\{ \gamma \sigma(X'_t)^{-1}dX'_t
+ \sqrt{1-\gamma^2}dU'_t  \right\},
\]
which is (\ref{stjerne-dif}) for $d=1$.

\vspace{3mm}

\subsection{Implementation}
\label{imp}

A main problem when implementing the proposed algorithm is to decide
whether or not coupling has happened, i.e.\ whether the two processes
have hit each other. In the following we discuss criteria for deciding whether
the two processes have met.

A criterion for deciding whether the two processes have met can be
found in the following manner. Only the case where ${\cal O} = I$ 
will be considered. We must determine whether the 
diffusions $X_t = X^*_{T-t}$ and $X'_t$ have coupled in a time interval 
$[\delta i,\delta (i+1)]$ (given that they did not meet before time 
$\delta i$). We will make the simplifying assumption that the diffusion 
$X$ develops according to (\ref{basicsde}), i.e.\ we ignore the influence 
of the fact that we have conditioned on $X_T=b$. If $\delta$ is 
sufficiently small and if the two diffusions are sufficiently close, 
we can assume that the drift and diffusion coefficients are constant 
in the time interval $[\delta i, \delta (i+1)]$ and that the coefficients 
are equal for the two processes. With these approximations 
\[
X'_{\delta i +s} = X'_{\delta i} + \alpha s 
+ \sigma \{ I - (1-\gamma)\Pi(X_{\delta i},X'_{\delta i})\} B_s 
+ \sqrt{1-\gamma^2}u (X_{\delta i},X'_{\delta i}) B^1_s, 
\]
where $B$ is a $d$-dimensional standard Wiener process, and $B^1$ is a 
one-dimensional standard Wiener process independent of $B$. Thus
$\sigma^{-1} (X_{\delta i +s} -\alpha s)$ and $\sigma^{-1} (X'_{\delta i +s} 
-\alpha s)$ are Brownian motions (started at $\sigma^{-1} X_{\delta i}$ 
and $\sigma^{-1} X'_{\delta i}$, respectively). The projection of these 
Brownian motions onto the plane $L$ are identical, where 
\[
L = \left\{ u \, | \, (u - \mbox{\small $\frac12$} \sigma^{-1} (X_{\delta i} +
X'_{\delta i}))^T \sigma^{-1} (X_{\delta i} - X'_{\delta i}) = 0 \right\}.
\]
This is the plane orthogonal to $\sigma^{-1}(X_{\delta i} - X'_{\delta i})$ 
through the point $\frac12 (X_{\delta i}+X'_{\delta i})$.
Therefore, if their projections onto the subspace orthogonal to
$L$ (i.e.\ the subspace generated by $\sigma^{-1}(X_{\delta i} -
X'_{\delta i})$) have passed each other at time $\delta (i+1)$, the two
$d$-dimensional processes must by continuity have been equal at some
time-point in $[\delta i, \delta (i+1)]$. Hence coupling must have
occurred in $[\delta i, \delta (i+1)]$ if 
\be
\label{crit1}
(X_{\delta i}- X'_{\delta i})^T V^{-1} (X_{\delta (i+1)} - 
X'_{\delta (i+1)}) < 0.
\ee
Here we have used that $\sigma^{-1}(X_{\delta i} - X'_{\delta i})$ is
orthogonal to $L$, and that $(X_{\delta i}- X'_{\delta i})^T
V^{-1} (X_{\delta (i)} -  X'_{\delta (i)}) > 0$. 

For the reflection method ($\gamma = -1$) we can find an alternative 
criterion and say a bit more. Using the same approximation, we find that
\bean
X_{\delta i +s} &=& X_{\delta i} + \alpha s + \sigma B_s \\
X'_{\delta i +s} &=& X'_{\delta i} + \alpha s 
+ \sigma H (X_{\delta i},X'_{\delta i}) B_s.
\eean
Thus $\sigma^{-1} (X_{\delta i +s} -\alpha s)$ is a Brownian motion started
at $\sigma^{-1} X_{\delta i}$, while $\sigma^{-1} (X'_{\delta i +s} 
-\alpha s)$ is the reflection of this Brownian motion in the plane $L$.
This follows from the fact that the
matrix $H (X_{\delta i},X'_{\delta i})$ is reflection in the plane
orthogonal to the vector $\sigma^{-1}(X_{\delta i} - X'_{\delta i})$. 
Note that $X_{\delta i +s} = X'_{\delta i +s}$ if and only if
$\sigma^{-1} (X_{\delta i +s} -\alpha s) = \sigma^{-1} (X'_{\delta i +s} 
-\alpha s)$, but as the latter two are each others reflection in the
plane $L$, they must meet in this plane if they intersect. Thus if
$\sigma^{-1} (X_{\delta i +s} -\alpha s)$ has crossed the plane $L$ in
the time interval $[\delta i,\delta (i+1)]$, then $X$ and $X'$ must
have intersected in this time interval. This is certainly the case if
\[
\{X_{\delta (i+1)} -\alpha \delta - \mbox{\small $\frac12$} 
(X_{\delta i} + X'_{\delta i})\}^T V^{-1} (X_{\delta i} -
X'_{\delta i}) \leq 0,
\]
where $V = \sigma \sigma^T$  or equivalently if
\be
\label{coupling}
X_{\delta (i+1)}^T V^{-1} (X_{\delta i} - X'_{\delta i})
\leq
\{\delta \alpha + \mbox{\small $\frac12$} (X_{\delta i} + 
X'_{\delta i})\}^T V^{-1} (X_{\delta i} - X'_{\delta i}).
\ee
In our simulation algorithm we can therefore assume
that coupling happens in the time interval $[\delta i,\delta (i+1)]$
if
\[
(Y^*_{\delta (N-i-1)})^T V^{-1} (Y^*_{\delta (N-i)} - Y'_{\delta i})
\leq
\{\delta \alpha (Y^*_{\delta (N-i)}) + \mbox{\small $\frac12$} 
(Y^*_{\delta (N-i)} 
+  Y'_{\delta i})\}^T V^{-1} (Y^*_{\delta (N-i)} - Y'_{\delta i}).
\]

Similar considerations can be used to estimate the probability that
coupling occurs in the time interval $[\delta i,\delta (i+1)]$. Since
(\ref{coupling}) implies coupling in $[\delta i,\delta (i+1)]$, the
probability of this event conditional on $X_{\delta i}$ and
$X'_{\delta i}$ is larger than the conditional probability that
\[
(X_{\delta (i+1)} - X_{\delta i} - \delta \alpha )^T V^{-1} 
(X_{\delta i} - X'_{\delta i})
\leq
- \mbox{\small $\frac12$} (X_{\delta i} - 
X'_{\delta i})^T V^{-1} (X_{\delta i} - X'_{\delta i}).
\]
Since (conditional on $X_{\delta i}$ and $X'_{\delta i}$) 
\[
(X_{\delta (i+1)} - X_{\delta i} - \delta \alpha )^T V^{-1} 
(X_{\delta i} - X'_{\delta i}) \sim \mbox{N}(0,\delta \, \omega^2),
\]
where $\omega^2 = (X_{\delta i} - X'_{\delta i}))^T V^{-1} 
(X_{\delta i} - X'_{\delta i})$, the conditional probability of
coupling in $[\delta i,\delta (i+1)]$ is larger than $\Phi (-\frac12
\omega/\sqrt{\delta})$, where $\Phi$ is the distribution function of
the standard normal distribution. 

In the MCMC algorithms to simulate exact diffusion bridges, it must be
determined whether the approximate bridge, $Z$, in a certain time step
has been hit by the associated $p^*_T(b)$-diffusion, $\tilde X$. Also
to determine whether this has happened the two criteria just derived
can be used. The reason is that the $p^*_T(b)$-diffusion is related to
$Z$ exactly as $X'$ is related to $X$. Therefore the same
approximations and calculations can be made.

\vspace{10mm}

\sect{Simulation study}
\label{simul}

In this section we test our simulation algorithm by applying it to 
the 2-dimensional Ornstein-Uhlenbeck bridge in the time interval
$[0,1]$, which can easily be simulated
exactly by another method, and for which the marginal distribution is known 
explicitly. We can therefore compare our method to exact results for the 
Ornstein-Uhlenbeck bridge. We do this by simulating a large number 
of $(a,b,1)-$bridges for selected values of $a$ and $b$. Then we compare the 
distribution of the simulated bridge to the exact distribution at 
time 0.5 (where the bridge-effect is strongest; see \cite{bladtsorensen}). The marginal distributions
are compared by Q-Q-plots, while the association between the coordinates
is evaluated by comparing the empirical copulas.

\vspace{3mm}

\subsection{The Ornstein--Uhlenbeck bridge}

In this subsection we give results on the $d$-dimensional 
Ornstein--Uhlenbeck process, which is given by
\be
\label{OU}
dX_t = -B(X_t - A)dt + \sigma dW_t,
\ee
where $A \in \R^d$, while $B$ and $\sigma$ are $d \times d$-matrices. It is
well-known that if all (complex) eigenvalues of $B$ have positive real
parts, then $X$ is ergodic with invariant density $N_d(A,\Gamma)$, the
$d$-dimensional normal distribution with mean $A$ and covariance matrix
\[
\Gamma = \int_0^\infty e^{-sB}Ve^{-sB^T}ds,
\]
where $V = \sigma \sigma^T$. The covariance matrix $\Gamma$ is the
unique symmetric solution to the equation
\be
\label{Gamma}
B \Gamma + \Gamma B^T = V.
\ee

The time-reversed version (\ref{reversesde}) 
of (\ref{OU}) has drift
\[
(B- V \Gamma^{-1})(X_t - A) = -\Gamma B^T \Gamma^{-1}(X_t - A),
\]
where we have used (\ref{Gamma}). It is thus an Ornstein--Uhlenbeck
process with a different drift matrix. The process (\ref{OU}) is
time-reversible if and only if $\Gamma B^T \Gamma^{-1} = B$ or 
$\Gamma B^T = B \Gamma$, i.e.\ if and only if the matrix $B \Gamma$ is
symmetric. By inserting this in (\ref{Gamma}) we obtain another
criterion for time-reversibility: $\Gamma = \frac12 B^{-1} V$. Thus
$B^{-1} V$ is necessarily symmetric. 

Now assume we have a general Ornstein--Uhlenbeck for which $B^{-1}V$ is
symmetric. Then $\Gamma = \frac12  B^{-1}V$ is a symmetric solution to
(\ref{Gamma}). Hence the process is time-reversible and the covariance
matrix of the invariant density function is $\Gamma$. We summarize in
the following lemma.

\begin{lem}\label{lemma:OU1}
An ergodic Ornstein--Uhlenbeck process (\ref{OU}) is time-reversible if and
only if $B^{-1}V$ is symmetric. If $B^{-1}V$ is symmetric, the
invariant distribution is the $N_d(A,\frac12 B^{-1} V)$-distribution.
\end{lem}

Because the Ornstein--Uhlenbeck process is Gaussian, the
Ornstein--Uhlenbeck bridge can be simulated by the method in the following
lemma, which also gives an expression for the marginal distribution. For 
simplicity, we assume that $A=0$.
 
 \begin{lem}\label{lemma:michael}
Generate $X_{t_0},X_{t_1}, \ldots X_{t_n},X_{t_{n+1}}$,
where $0=t_0 < t_1 < \cdots < t_n < t_{n+1}$, by $X_0 = x_0$ and
\[
X_{t_i} = e^{-B(t_i - t_{i-1})}X_{t_{i-1}} + W_i, \ \ i =1, \ldots, n+1
\]
where the $W_i$s are independent and
\[
W_i \sim N\left( 0, \Gamma_{t_i - t_{i-1}} \right)
\]
with
\[
\Gamma_t = \int_0^t e^{-sB}Ve^{-sB^T}ds.
\]
Define
\[
Z_{t_i} = X_{t_i} + e^{-B(t_{n+1} - t_{i})} \Gamma_{t_i} 
\Gamma_{t_{n+1}}^{-1}(x - X_{t_{n+1}}), \ \ i=0, \ldots, n+1.
\]
Then $(Z_{t_0}, Z_{t_1}, \ldots , Z_{t_n}, Z_{t_{n+1}})$ is
distributed like an Ornstein--Uhlenbeck bridge with $Z_{t_0}=x_0$ and
$Z_{t_{n+1} = x}$.

For an Ornstein--Uhlenbeck bridge $(Z_t)$ in $[0,1]$ with $Z_0=x_0$ and
$Z_1=x$, the distribution of $Z_t$ is a $d$-dimensional normal
distribution with expectation
\[
e^{-Bt}x_0+e^{-B(1-t)}\Gamma_t \Gamma_1^{-1}(x-e^{-B}x_0)
\]
and covariance matrix
\[
\Gamma_t - e^{-B(1-t)}\Gamma_t \Gamma_1^{-1}\Gamma_t e^{-B(1-t)}.
\]
\end{lem}

The lemma follows by standard arguments from the fact that for the
Ornstein--Uhlenbeck process the conditional distribution of $X_t$ given $X_s$ ($s <
t$) is the $N_d(e^{-B(t-s)}, \Gamma_{t-s})$-distribution.

\vspace{3mm}

\subsection{Simulations}

We simulated bridges for the Ornstein--Uhlenbeck process with
$A=0$, $\sigma = I$ and
\[ 
B=\left\{ 
\begin{array}{ll}
1.5 & 1 \\
1 & 1.5 
\end{array}
\right\}.
\]
This process is ergodic and time-reversible with stationary distribution
$N_2(0,\Gamma)$, where
\[ 
\Gamma=\mbox{\small $\frac12$} B^{-1} = \left\{ 
\begin{array}{rr}
0.6 & -0.4 \\ 
-0.4 & 0.6 
\end{array}
\right\}.
\]
The diffusion bridges were simulated over the time interval $[0,1]$
using the Euler scheme with discretization level $N = 50$ (step
size $\delta=0.02$). For the Ornstein--Uhlenbeck process the Euler
scheme is equal to the Milstein scheme. The two diffusions $X$ 
and $X'$ were assumed to have intersected in the time interval 
$[\delta i, \delta (i+1)]$ if both $|X_{\delta i} - X'_{\delta i}| \leq 0.05$ and 
(\ref{crit1}) were satisfied.

\begin{figure}
\begin{center}
\vspace{-2cm}
\includegraphics[width=5.6cm]{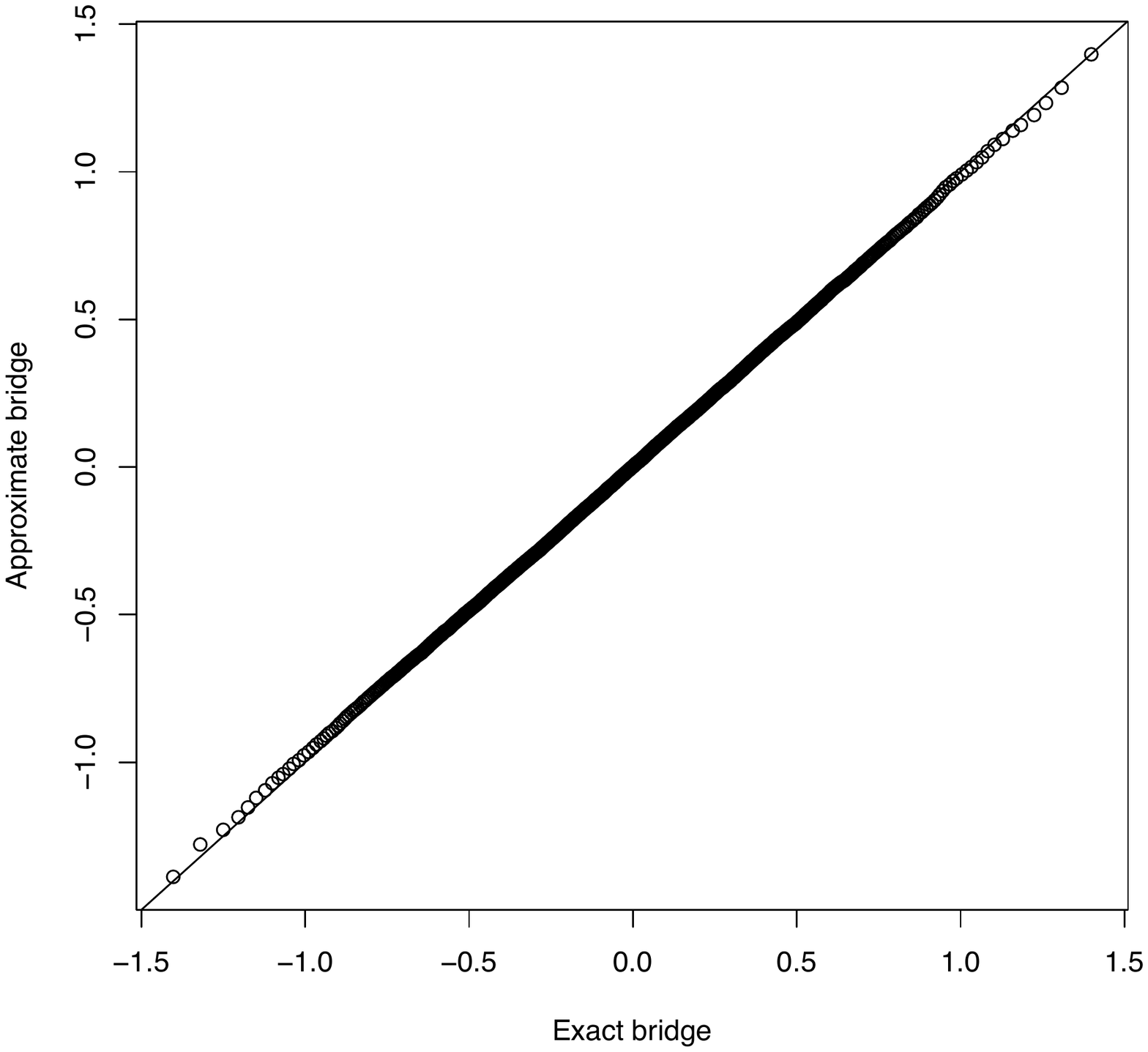} 
\hspace{-5mm}
\includegraphics[width=5.6cm]{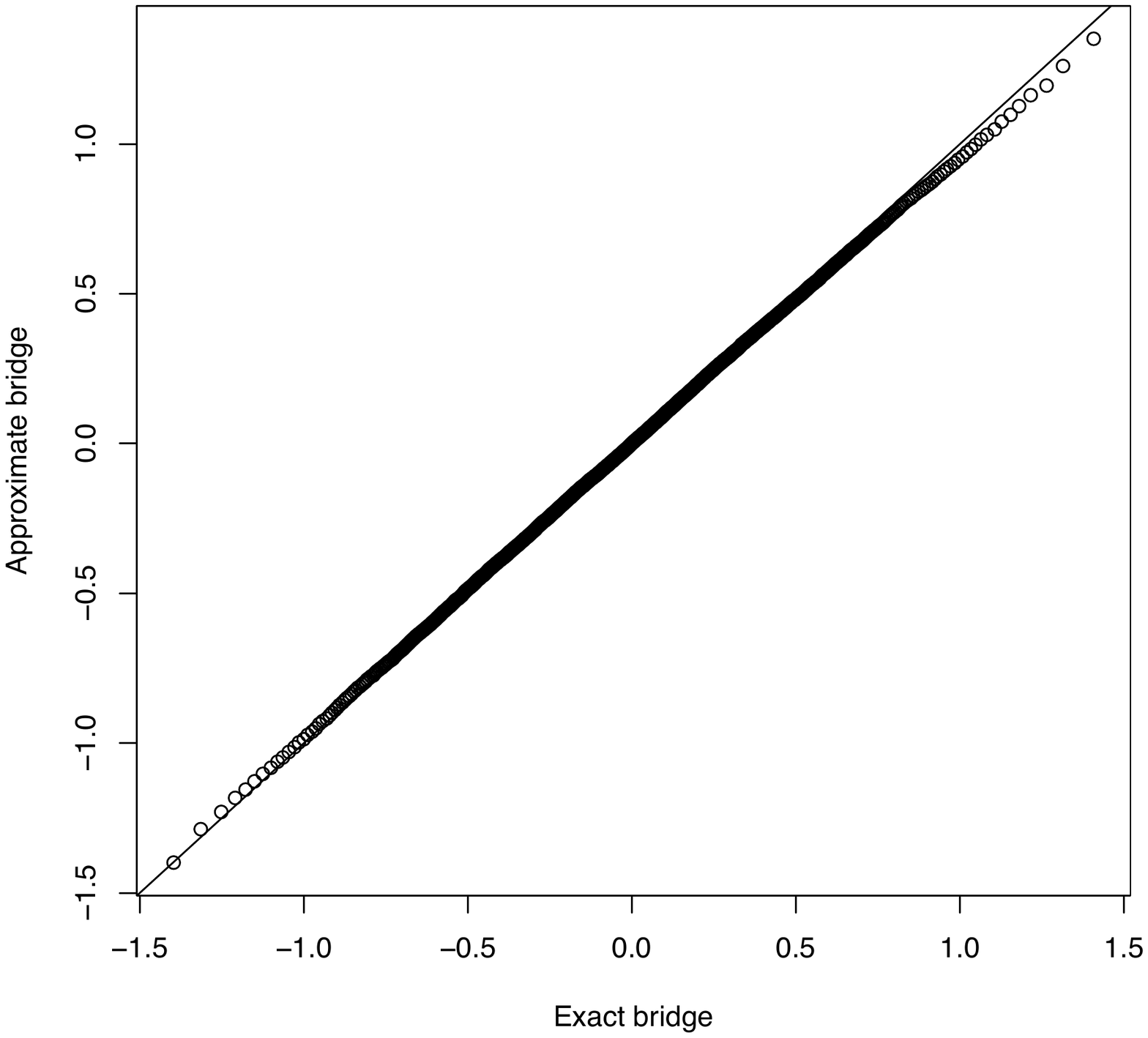}
\hspace{-5mm}
\includegraphics[width=5.6cm]{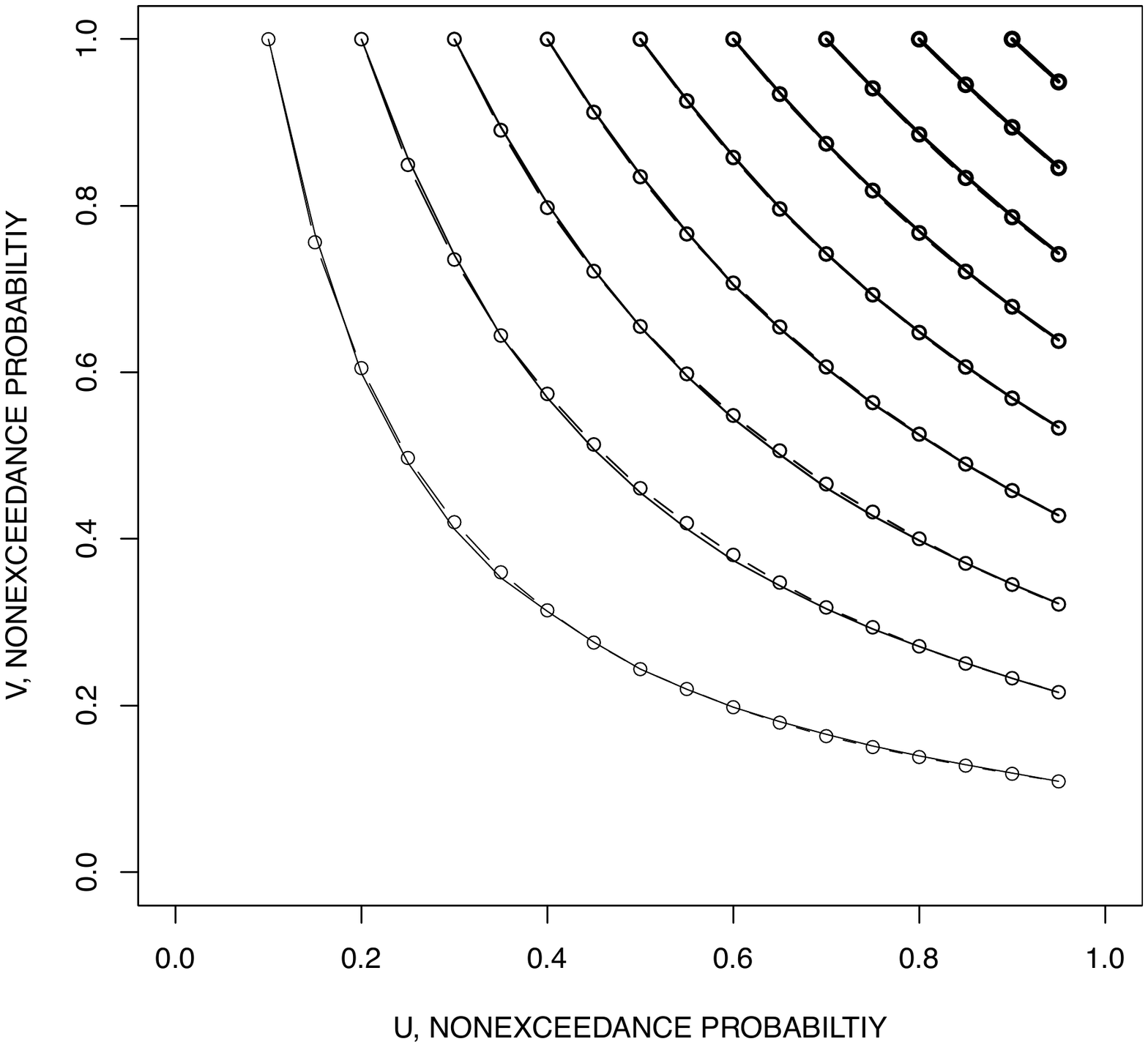}  
\vspace{-1.5cm}
\caption{\label{simulation1} 
Q-Q plots comparing the marginal distributions at time 0.5 for 
50.000 two-dimensional Ornstein--Uhlenbeck bridges from $(0,0)$ 
to $(0,0)$ simulated by the approximate method with $\gamma = -1$ 
(method of reflection) to the exact marginal distributions of 
Ornstein--Uhlenbeck bridges. Level curves of the empirical copula 
for the two 2-dimensional distribution at time 0.5 is compared to
those of the exact copula (full drawn curves).} 
\end{center}
\end{figure}

First we simulated 50.000 diffusion bridges from $(0,0)$ to $(0,0)$
using the approximate method presented in Section \ref{approxbridge} 
with $\gamma = -1$ (method of reflection). The computing time was 
6 seconds.
In Figure \ref{simulation1} the marginal distributions at time 0.5 for
the simulations are compared to the exact distributions given by 
Lemma \ref{lemma:OU1} by means of Q-Q plots. The dependence between 
the marginals (at time 0.5) in the simulations is compared to the exact 
dependence for an Ornstein--Uhlenbeck bridge by plotting the
level curves of the empirical and the exact copulas.
The two-dimensional
distribution at time 0.5 for the approximate bridge simulation method
is essentially equal to the distribution for the exact diffusion
bridge. Since the distribution of the approximate diffusion bridge
fits the distribution of the exact bridge, the MCMC method presented
in Section \ref{exactsimulation} does not improve the distribution of
the approximate bridge. Therefore the result is not plotted. 

\begin{figure}
\begin{center}
\vspace{-2cm}
\includegraphics[width=5.5cm]{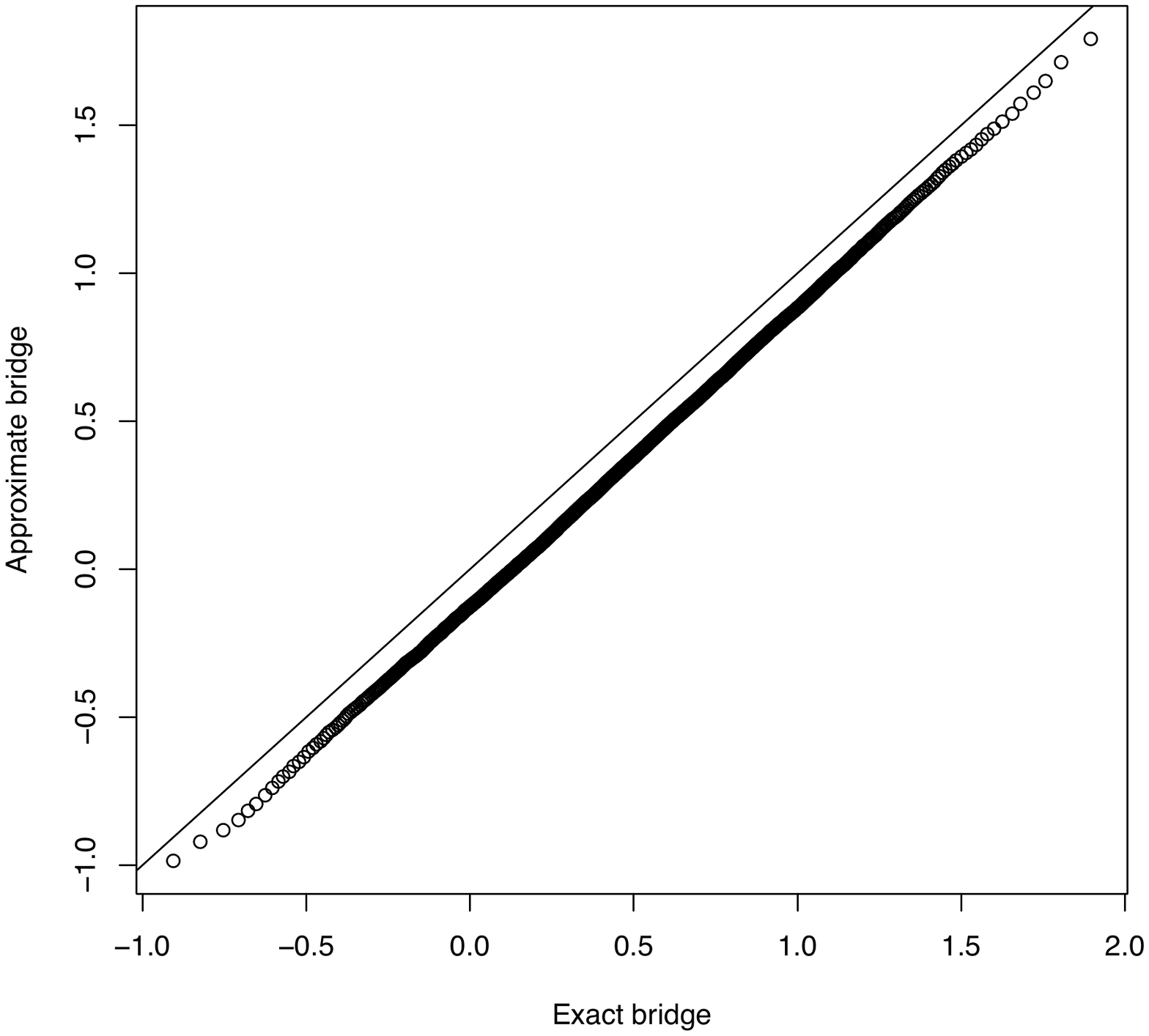} 
\hspace{-5mm}
\includegraphics[width=5.5cm]{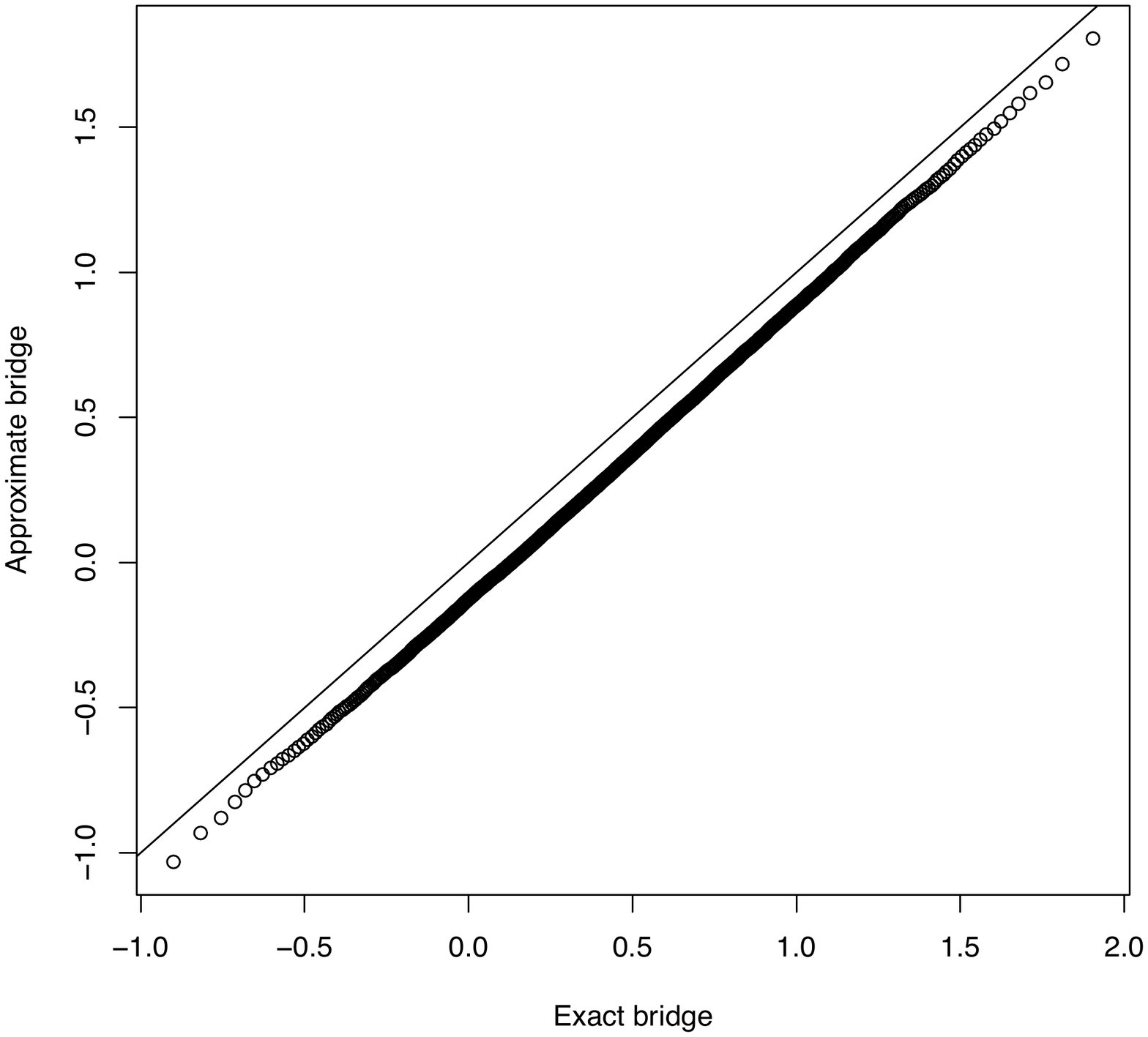}
\vspace{-5mm} 
\includegraphics[width=5.5cm]{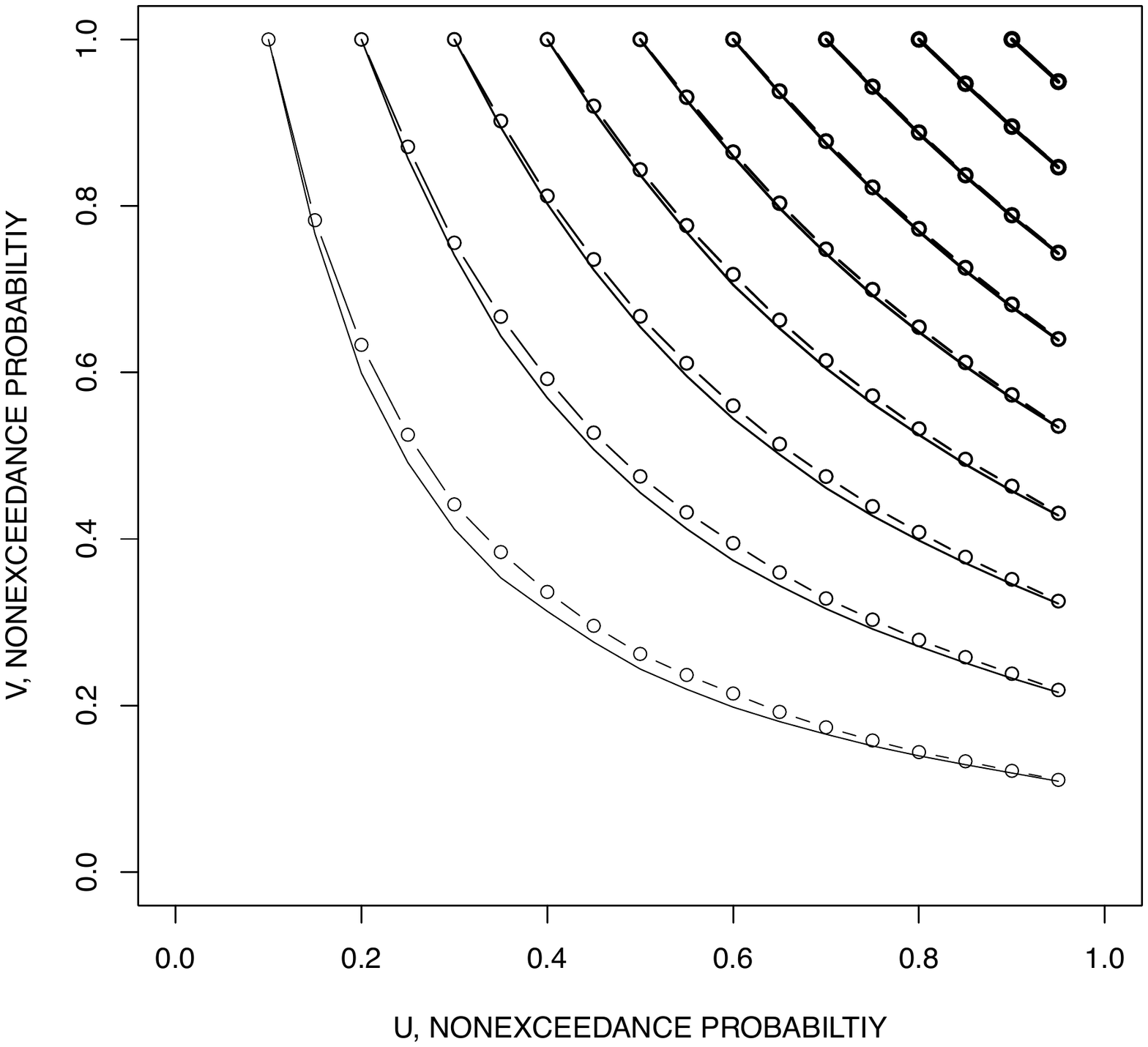} 
\vspace{-1.3cm}
\caption{\label{simulation2} 
Q-Q plots comparing the empirical marginal distributions at time 0.5 for 
50.000 simulated two-dimensional Ornstein--Uhlenbeck bridges from 
$(0.785,0.785)$
to $(1.091,1.091)$ by the approximate method with $\gamma = 0$ (method of 
projection) to the exact marginal distributions of Ornstein--Uhlenbeck 
bridges. Level curves of the empirical copula for the two 2-dimensional 
distribution at time 0.5 are compared to those of the exact copula 
(full drawn curves).} 
\end{center}
\end{figure}

\begin{figure}
\begin{center}
\vspace{-2.2cm}
\includegraphics[width=5.6cm]{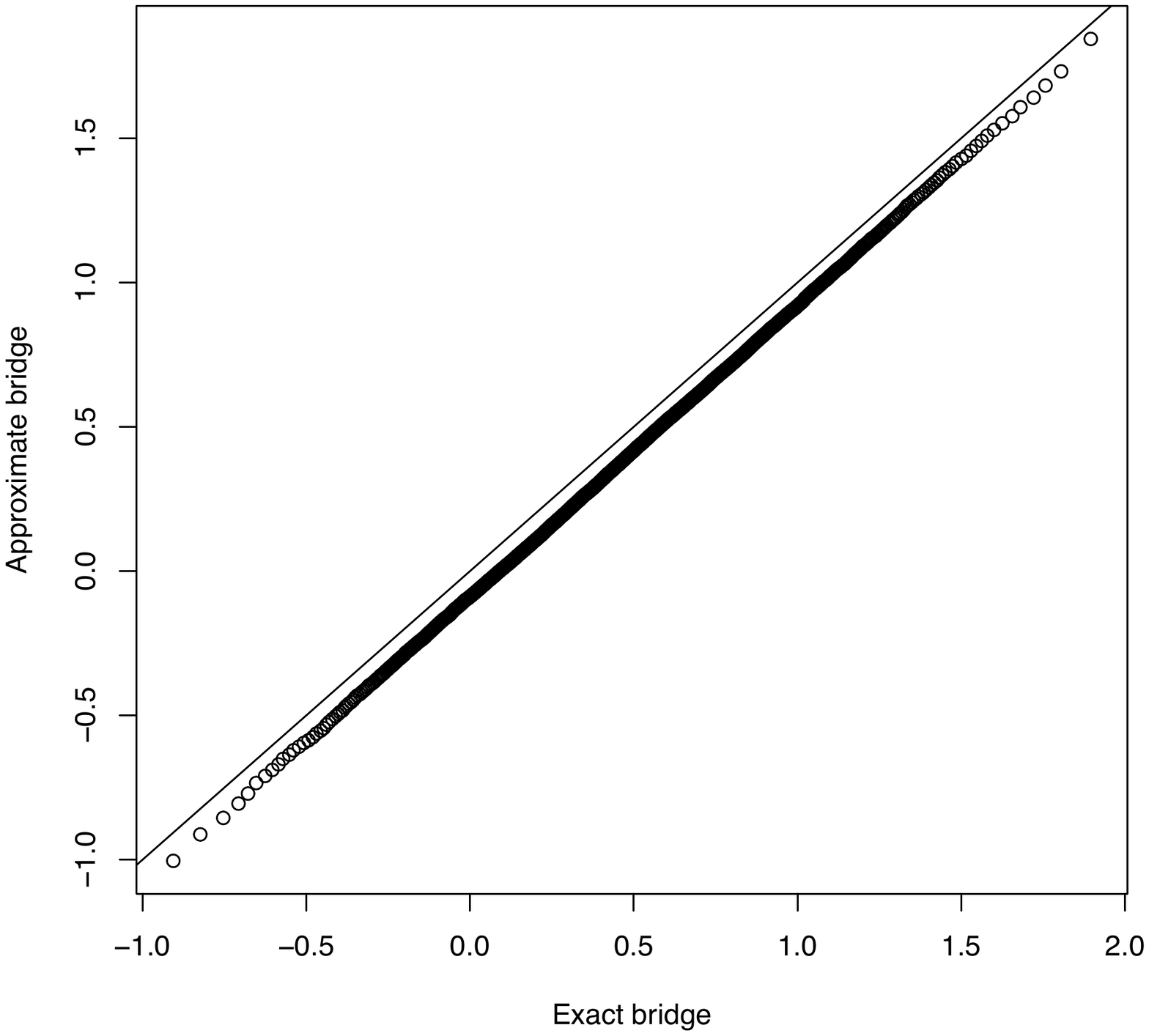} 
\hspace{-5mm}
\includegraphics[width=5.6cm]{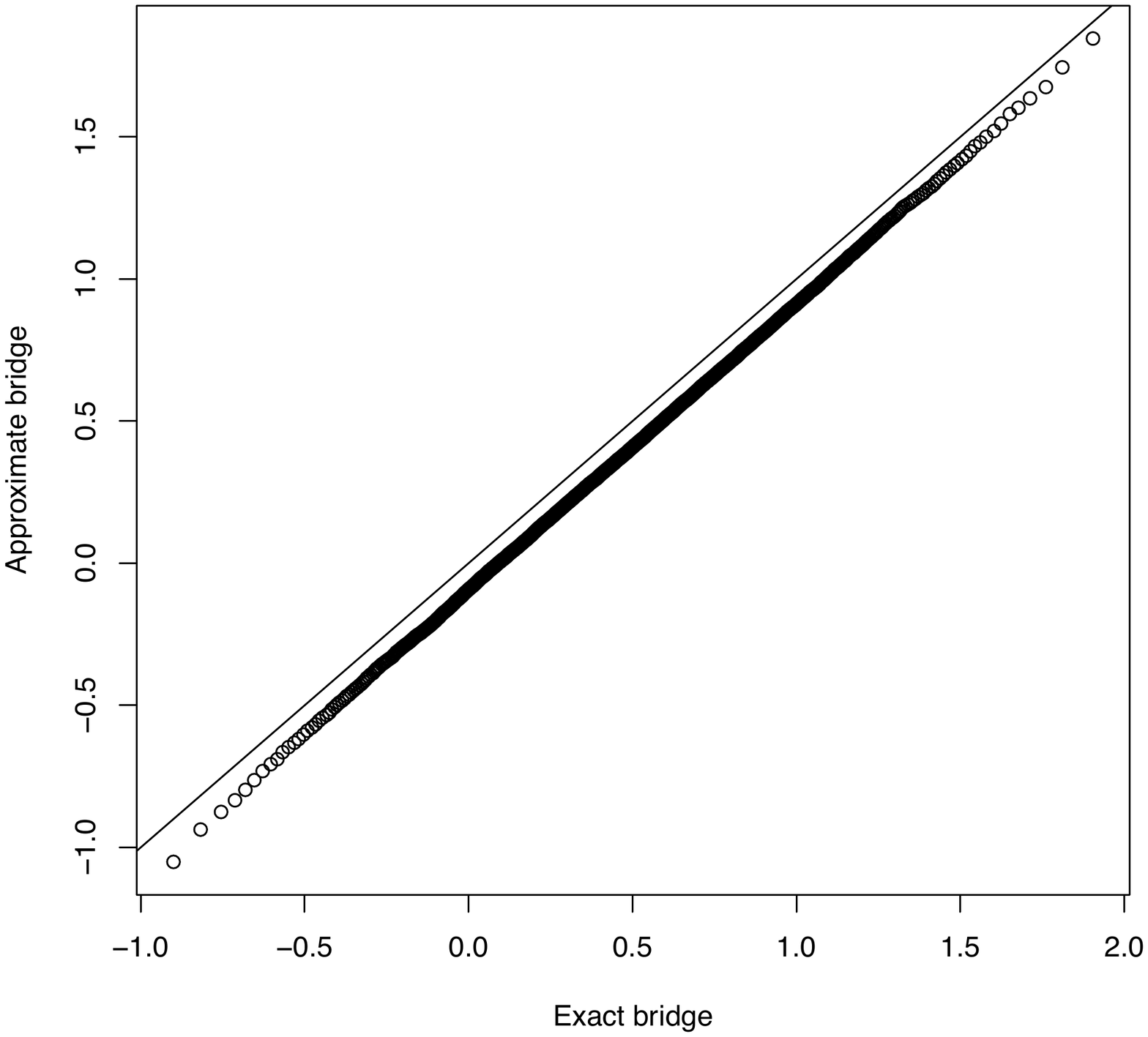}
\hspace{-5mm} 
\includegraphics[width=5.6cm]{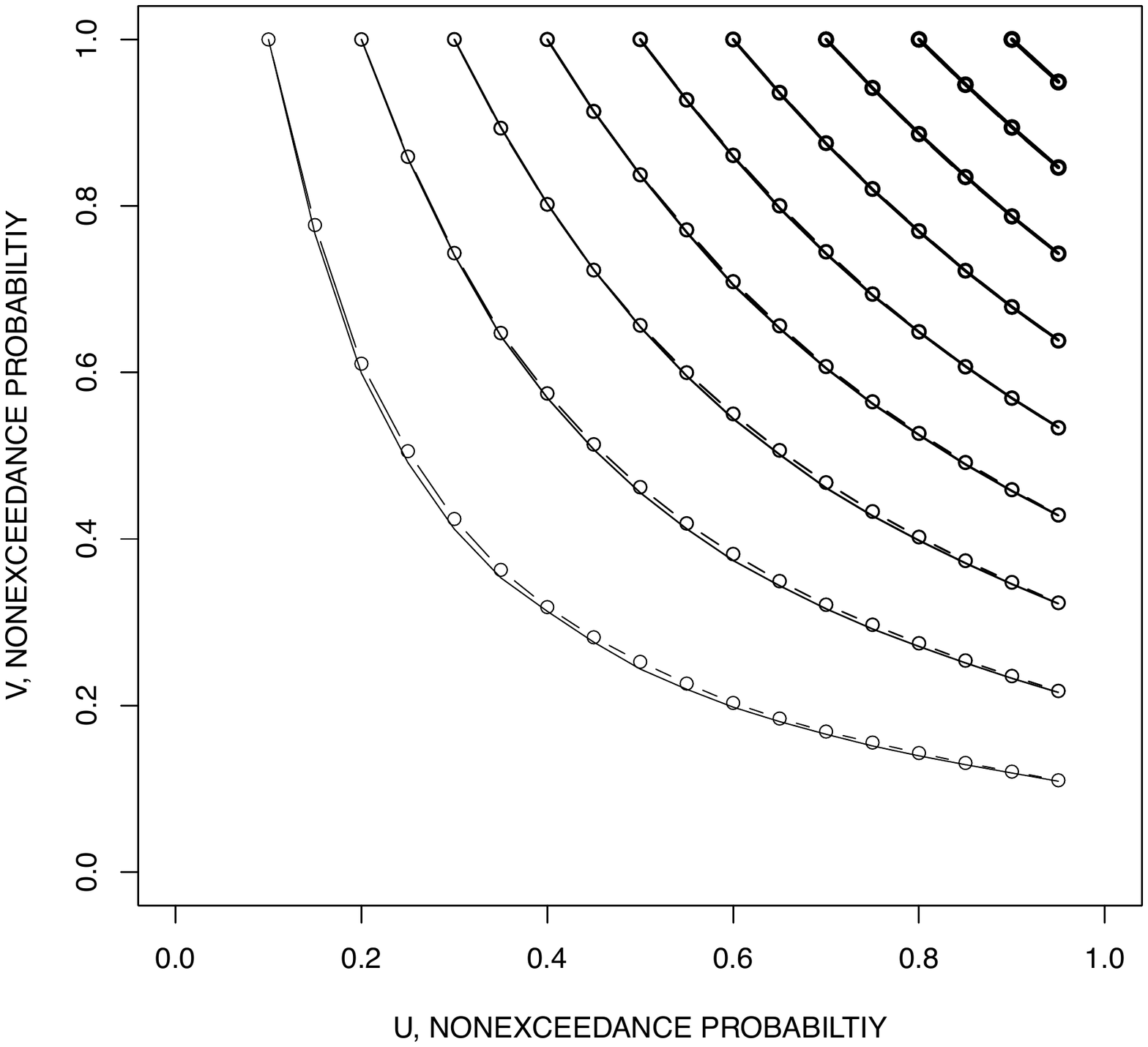} 
\vspace{-1.5cm}
\caption{\label{simulation3} 
Q-Q plots comparing the empirical marginal distributions at time 0.5 for 
50.000 simulated two-dimensional Ornstein--Uhlenbeck bridges from 
$(0.785,0.785)$
to $(1.091,1.091)$ by the approximate method with $\gamma = 0.5$ to
the exact marginal distributions of Ornstein--Uhlenbeck bridges. Level 
curves of the empirical copula for the two 2-dimensional distribution 
at time 0.5 are compared to those of the exact copula (full drawn 
curves).} 
\end{center}
\end{figure}

\begin{figure}
\begin{center}
\vspace{-2.5cm}
\includegraphics[width=5.5cm]{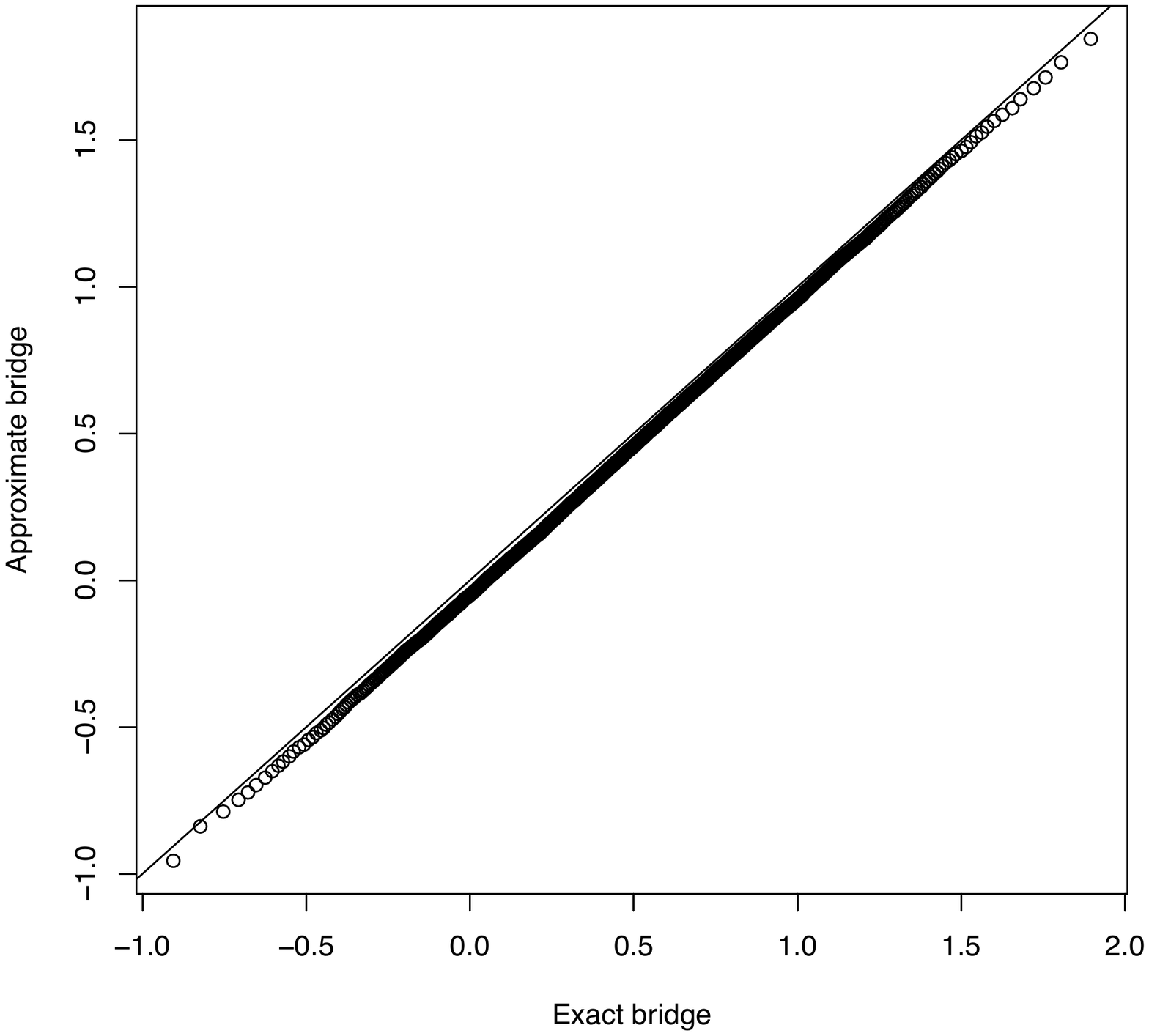} 
\hspace{-5mm}
\includegraphics[width=5.5cm]{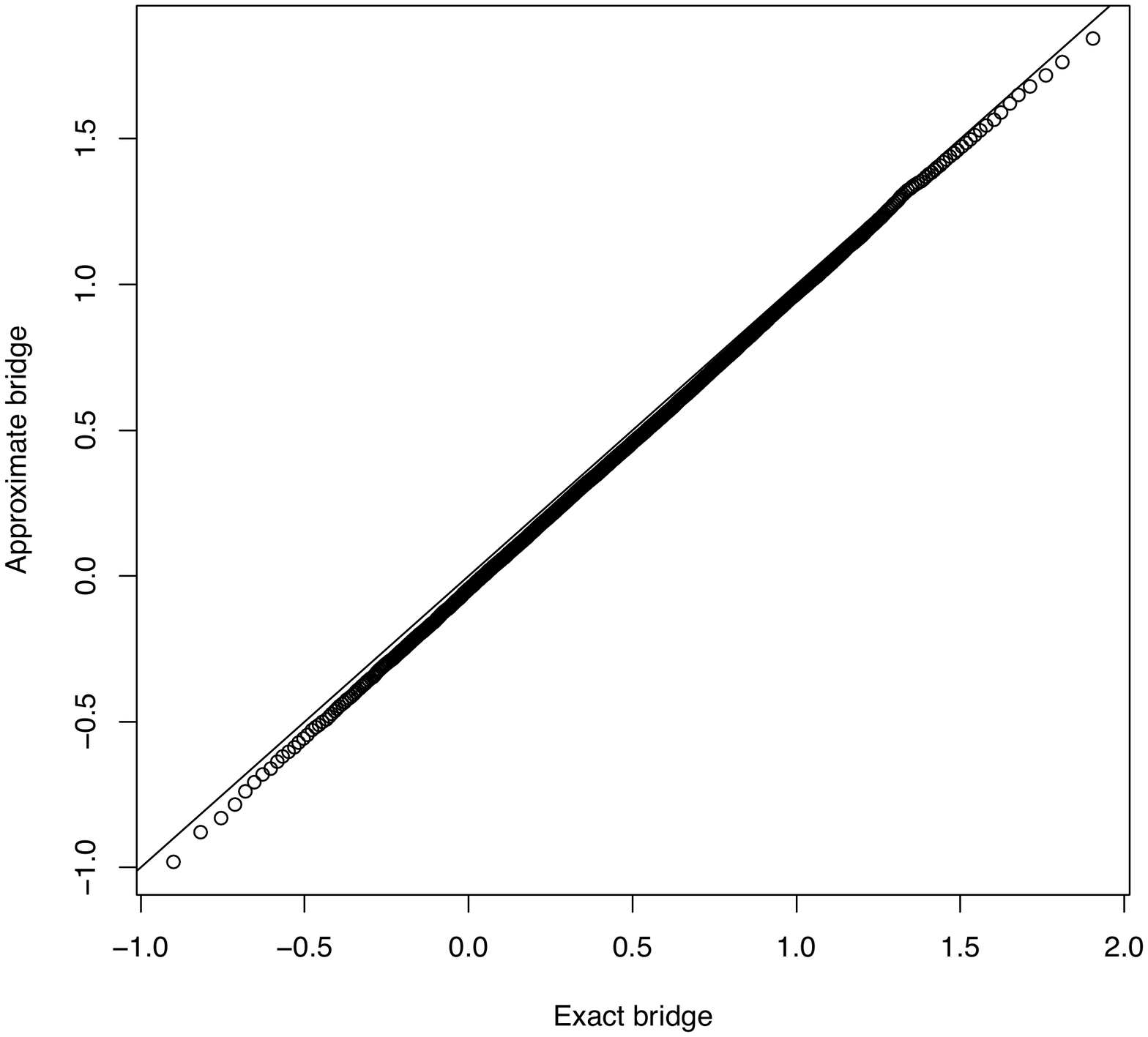}
\vspace{-5mm} 
\includegraphics[width=5.5cm]{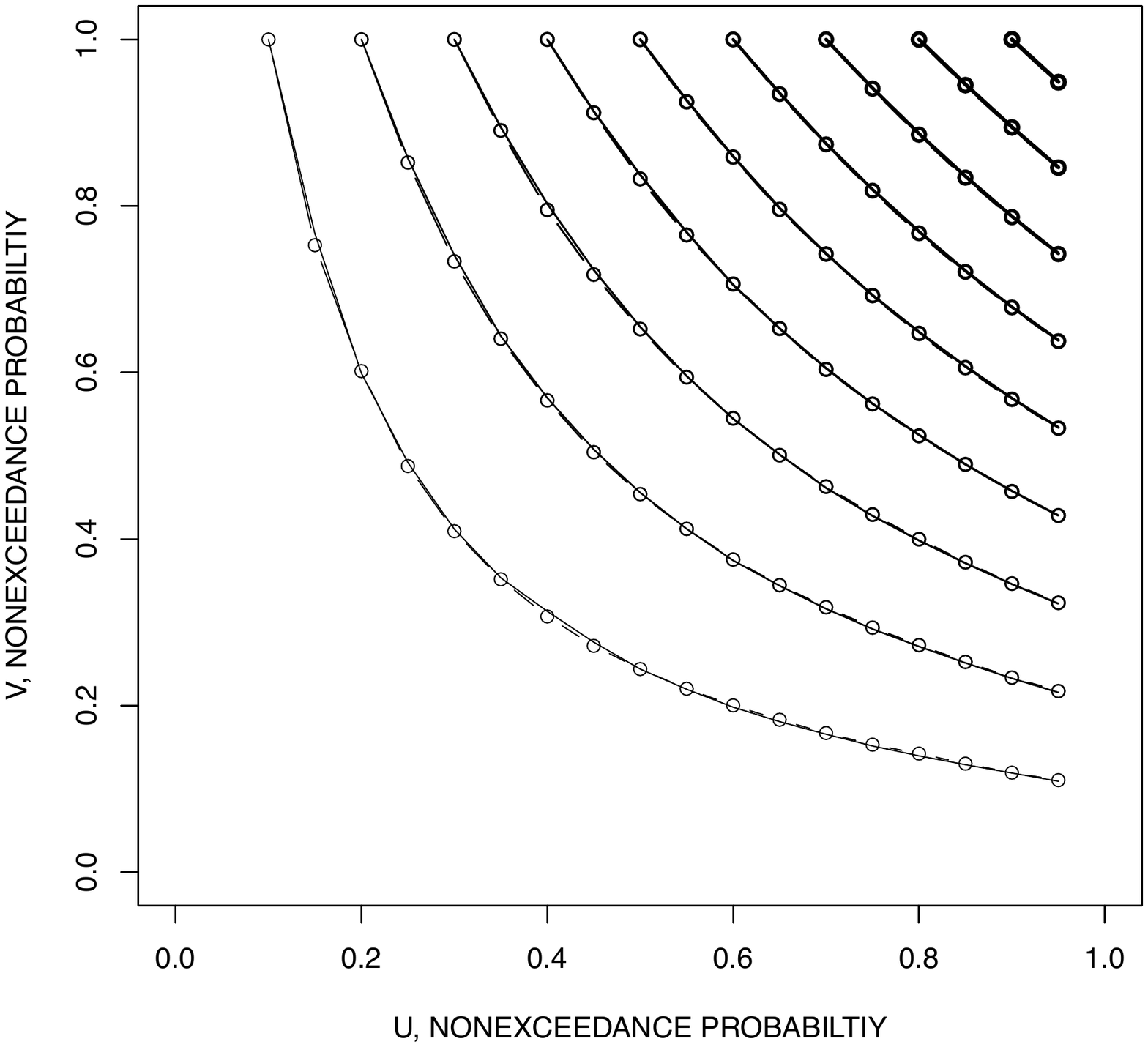} 
\vspace{-1.5cm}
\caption{\label{simulation4} 
Q-Q plots comparing the empirical marginal distributions at time 0.5 for 
50.000 simulated two-dimensional Ornstein--Uhlenbeck bridges from 
$(0.785,0.785)$
to $(1.091,1.091)$ by the approximate method with $\gamma = 0.9$ to
the exact marginal distributions of Ornstein--Uhlenbeck bridges. Level 
curves of the empirical copula for the two 2-dimensional distribution 
at time 0.5 are compared to those of the exact copula (full drawn curves).} 
\end{center}
\end{figure}

\begin{figure}
\begin{center}
\vspace{-2.2cm}
\includegraphics[width=5.5cm]{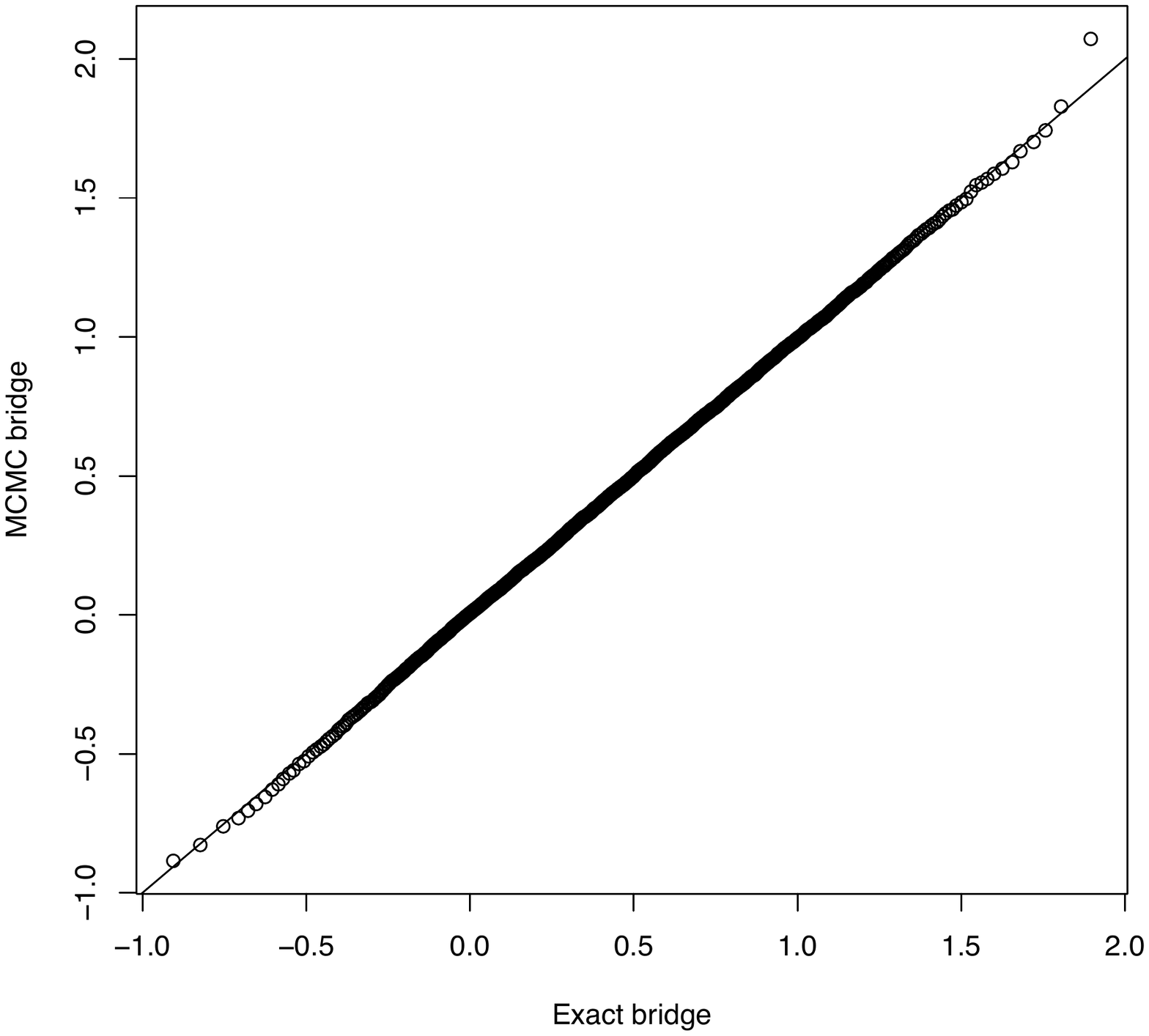} 
\hspace{-5mm}
\includegraphics[width=5.5cm]{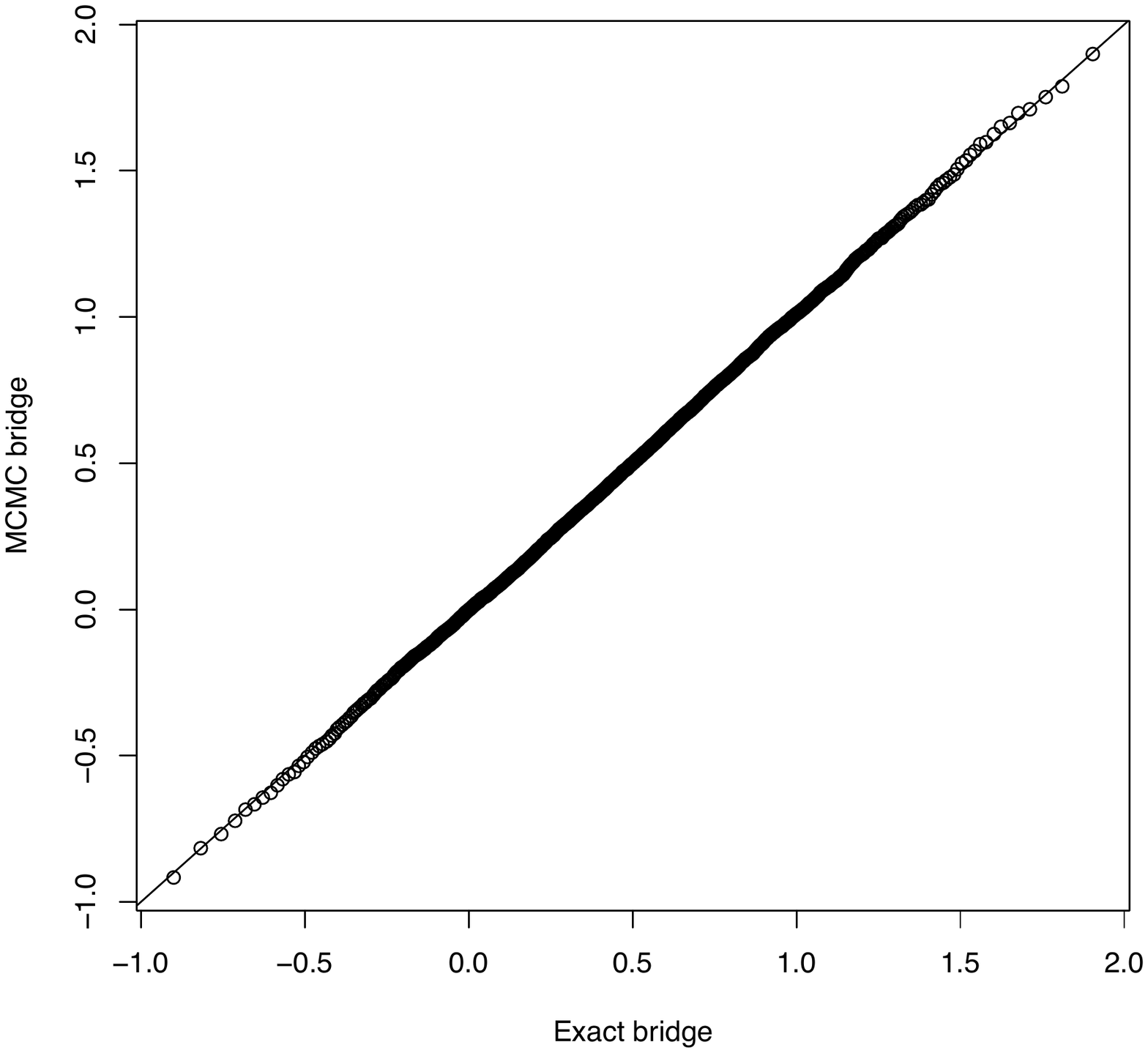}
\vspace{-5mm} 
\includegraphics[width=5.5cm]{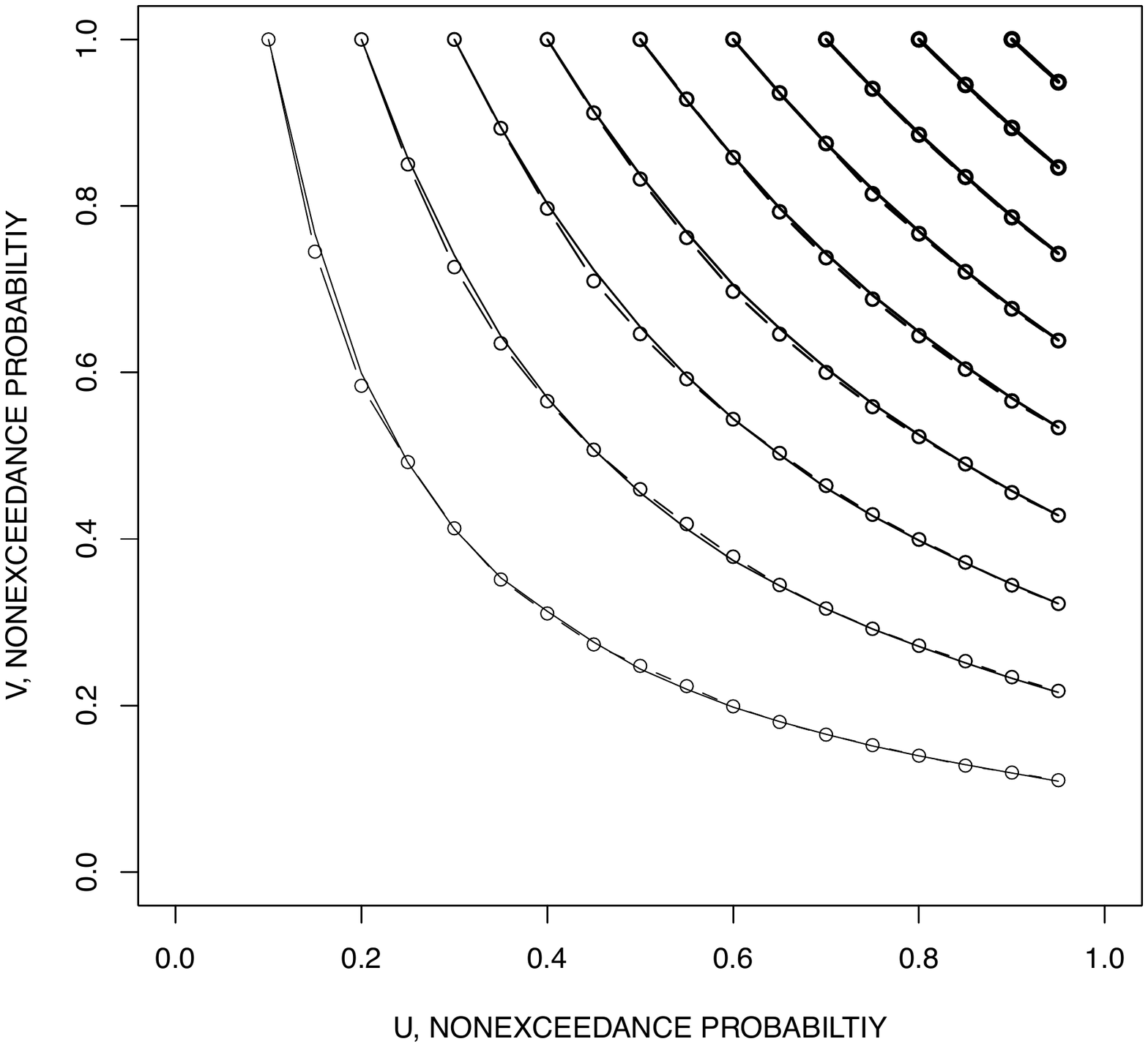} 
\vspace{-1.5cm}
\caption{\label{simulation5} 
Q-Q plots comparing the empirical marginal distributions at time 0.5 for 
two-dimensional Ornstein--Uhlenbeck bridges from $(0.785,0.785)$
to $(1.091,1.091)$ produced by 50.000 iterations of 
the pseudo-marginal Metropolis-Hastings algorithm with $\gamma = 0.5$ 
and $N=1$ to the exact marginal distributions of 
Ornstein--Uhlenbeck bridges. Level curves of the empirical copula for 
the two 2-dimensional distribution at time 0.5 are compared to those 
of the exact copula (full drawn curves).} 
\end{center}
\end{figure}

\begin{figure}
\begin{center}
\vspace{-2.3cm}
\includegraphics[width=5.5cm]{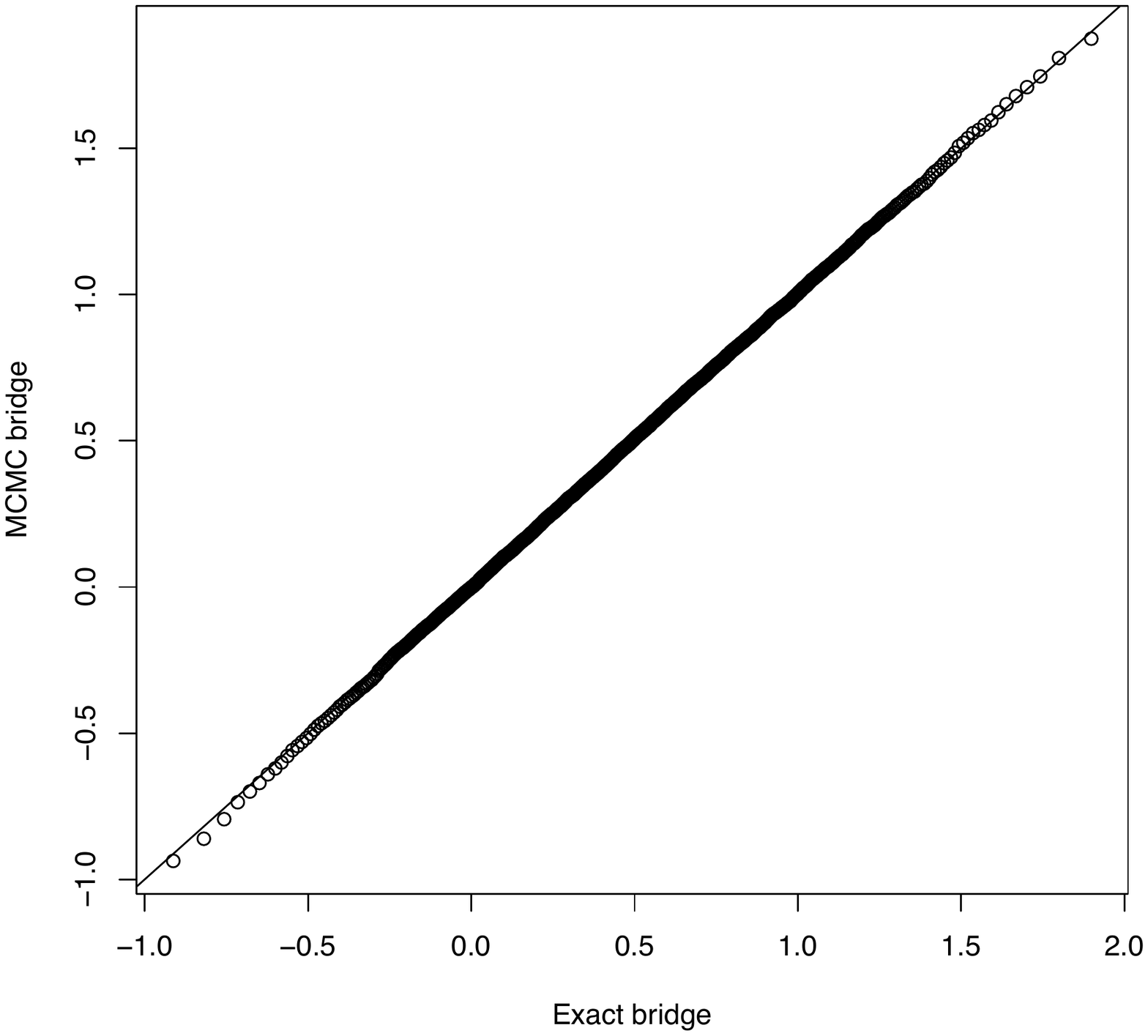} 
\hspace{-5mm}
\includegraphics[width=5.5cm]{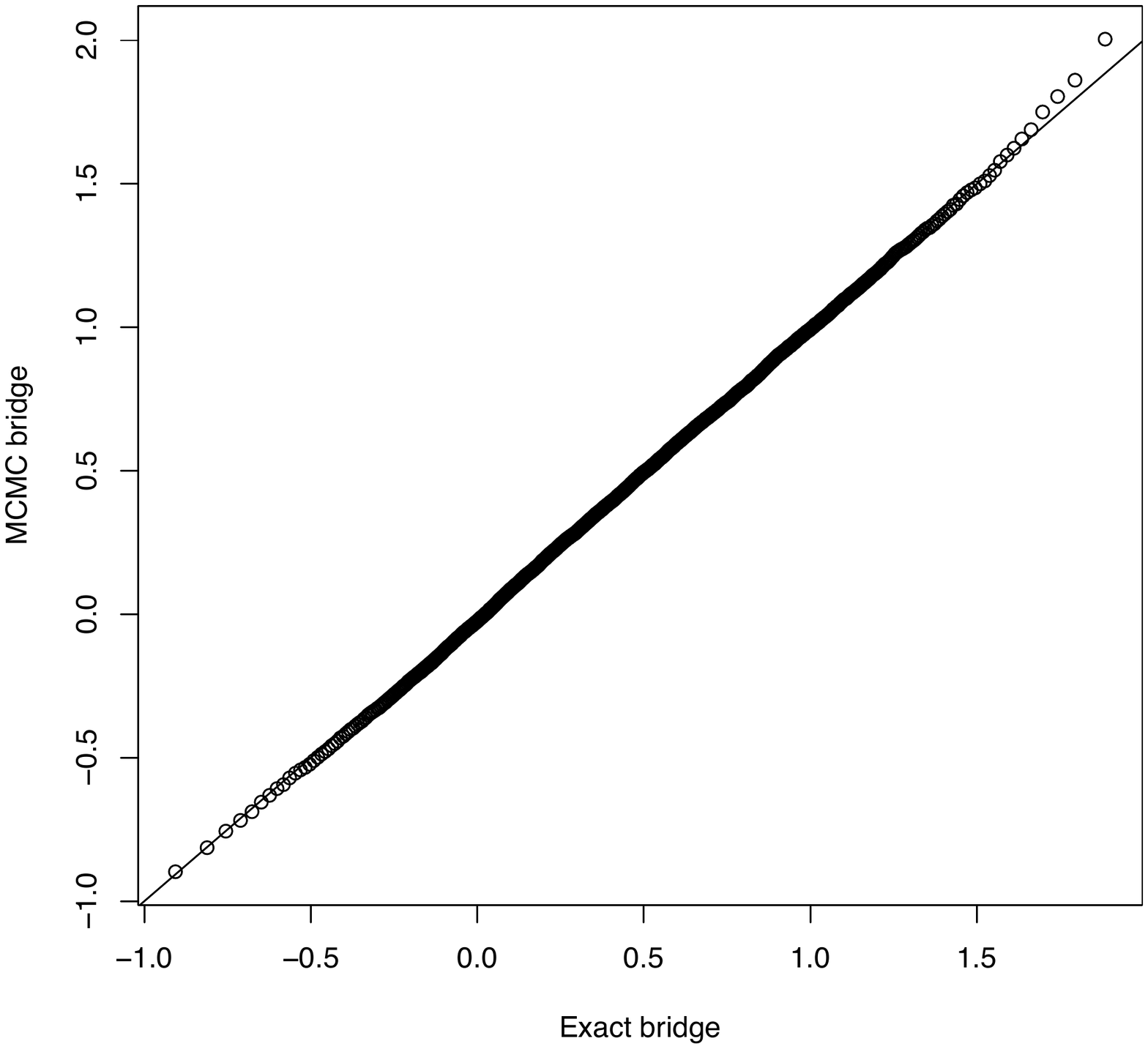}
\vspace{-5mm} 
\includegraphics[width=5.5cm]{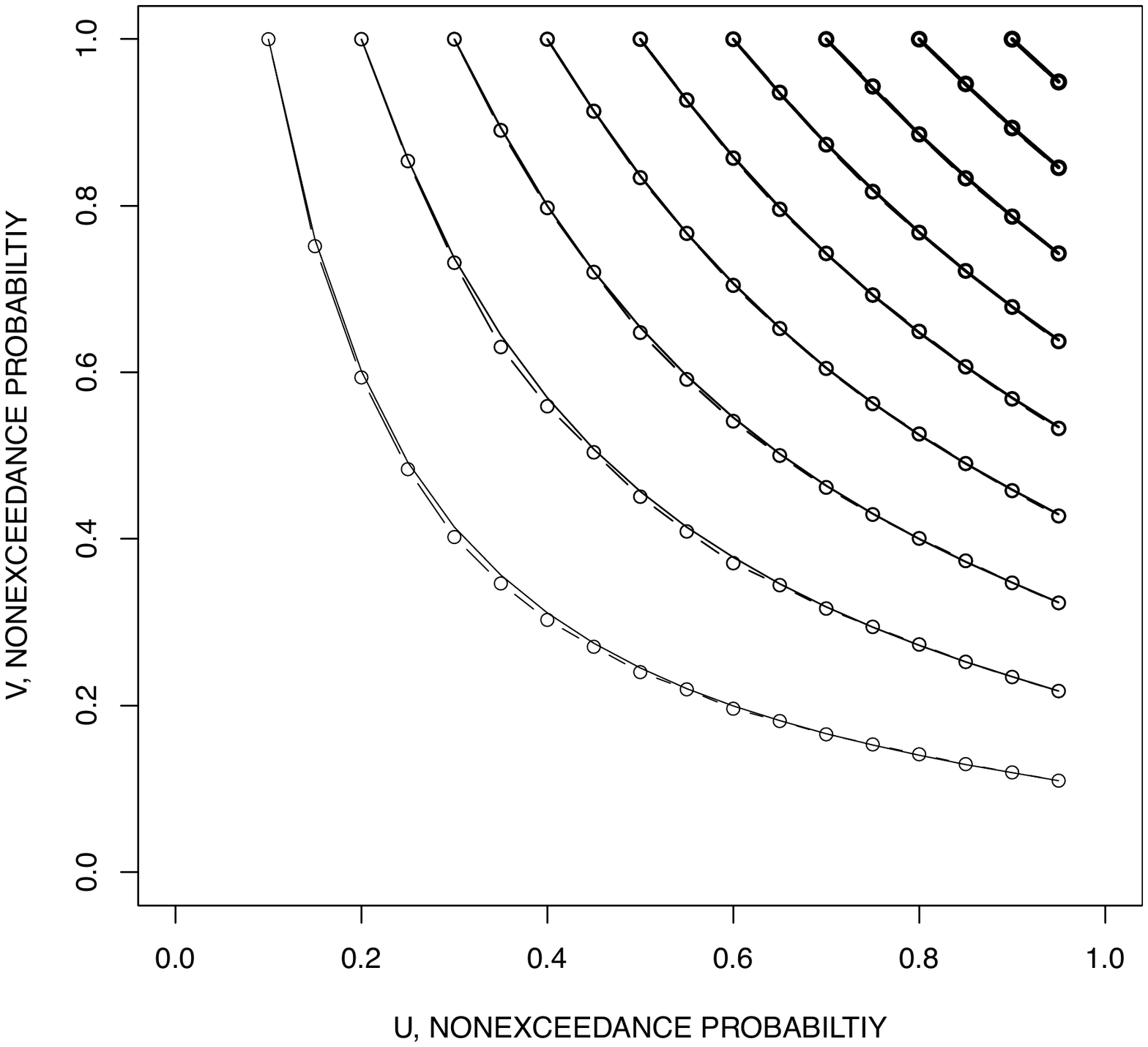} 
\vspace{-1.5cm}
\caption{\label{simulation6} 
Q-Q plots comparing the empirical marginal distributions at time 0.5 for 
two-dimensional Ornstein--Uhlenbeck bridges from $(0.785,0.785)$
to $(1.091,1.091)$ produced by 500.000 iterations of the alternative 
MCMC method with $\gamma = 0.5$ to the exact marginal distributions of 
Ornstein--Uhlenbeck bridges. Level curves of the empirical copula for the 
two 2-dimensional distribution at time 0.5 are compared to 
to those of the exact copula (full drawn curves).} 
\end{center}
\end{figure}

Similarly nice fits can be produced for diffusion bridges between points that, 
like $(0,0)$, do not have a small probability of being reached by the
diffusion. To test the method for a more unlikely bridge, diffusion
bridges from $(0.785,0.785)$ to $(1.091,1.091)$ were simulated. These
are bridges from the boundary of the $95.5\%$-ellipse of the
stationary distribution to the boundary of the $99.7\%$-ellipse. Such
bridges are rarely observed in data and are thus rarely needed for
simulated likelihood or Bayesian estimation. We simulated diffusion 
bridges using the approximate method in Section \ref{approxbridge} with
$\gamma = 0$, $\gamma = 0.5$ and $\gamma = 0.9$. For each of the 
$\gamma$-values we simulated 50.000 bridges. Marginal 
distributions and the copulas (at time 0.5) are compared to exact 
results (Lemma \ref{lemma:OU1}) in the Figures \ref{simulation2},
\ref{simulation3} and \ref{simulation4}.
The distributions of the approximate bridges do not exactly fit the
distribution of the  Ornstein-Uhlenbeck bridge for this unlikely
bridge, but the fit is rather good, and improves as $\gamma$
increases. For $\gamma$ very close to one (e.g. 0.95), the fit is essentially
exact. By using the pseudo-marginal Metropolis-Hastings algorithm in
Section \ref{exactsimulation} an essentially exact fit is obtained when
using any of approximate bridges with $\gamma = 0$, $\gamma = 0.5$ and
$\gamma = 0.9$. As an example Figure \ref{simulation5} compares the marginal 
distributions and the copula (at time 0.5) produced by the
pseudo-marginal Metropolis-Hastings algorithm where the proposal is
the approximate bridge with $\gamma = 0.5$. We ran 50.000 iterations of the 
algorithm with $N=1$. The computing time was 44 seconds. Generally,
the computing time for generating 50.000 bridges with the
pseudo-marginal Metropolis-Hastings algorithm varied from 40 to 120
seconds and increased with the value of $\gamma$.
Figure \ref{simulation6} compares the marginal distributions and the copula
(at time 0.5) produced by the alternative MCMC algorithm presented in
Section \ref{exactsimulation}. Again the proposal is the approximate 
bridge with $\gamma = 0.5$. The alternative MCMC algorithm has a much larger
rejection probability than the MH-algorithm, so to obtain plots of about 
the same quality as in Figure \ref{simulation5}, it was necessary to run 
500.000 iterations of the algorithm. This produced a sample with a large 
number of identical bridges, so it seems advisable to use only a subset of
the output from this algorithm, e.g.\ every 10th bridge. The computing 
time was 32 seconds. 

\begin{figure}
\begin{center}
\vspace{-2.2cm}
\includegraphics[width=5.5cm]{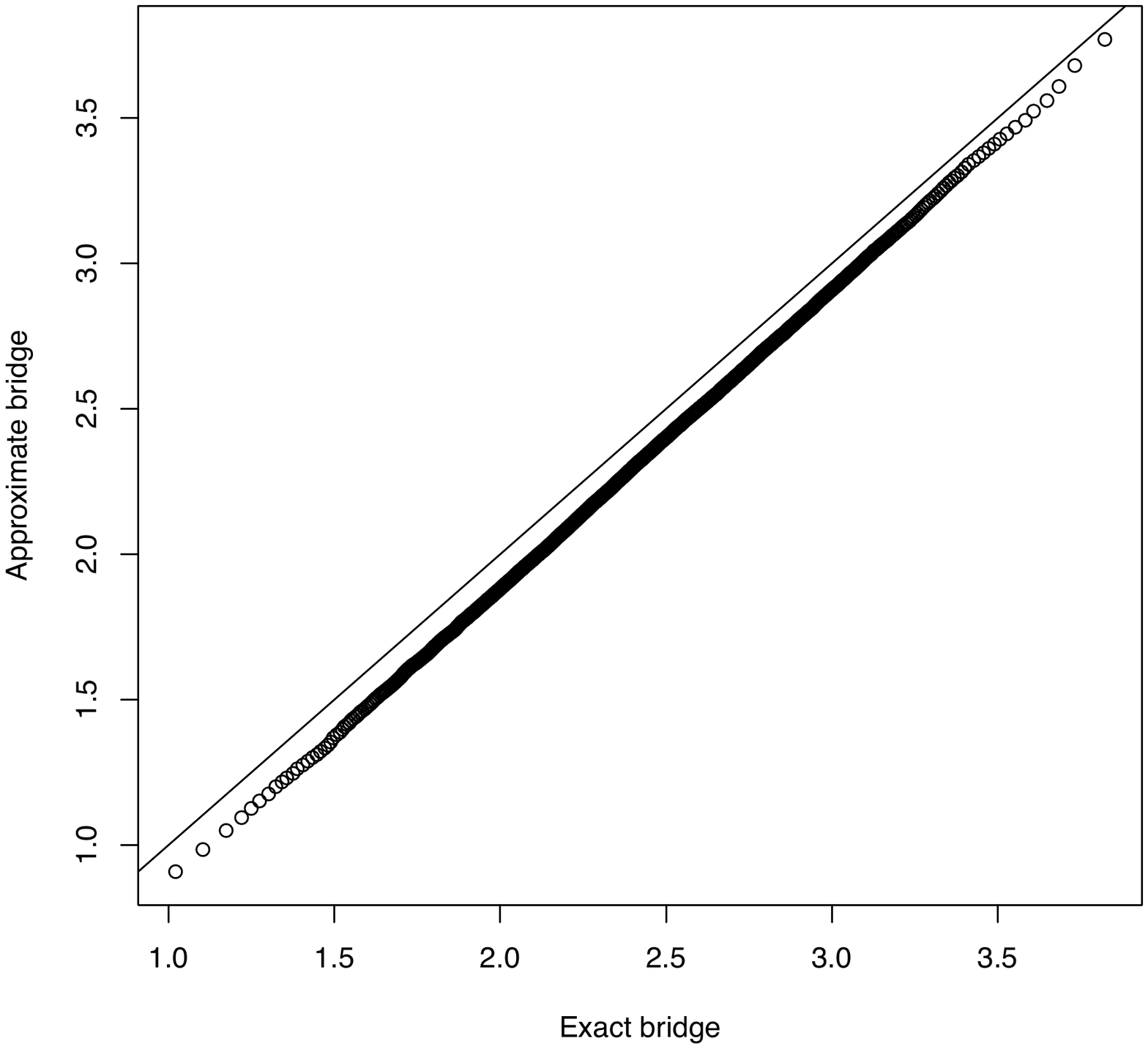} 
\hspace{-5mm}
\includegraphics[width=5.5cm]{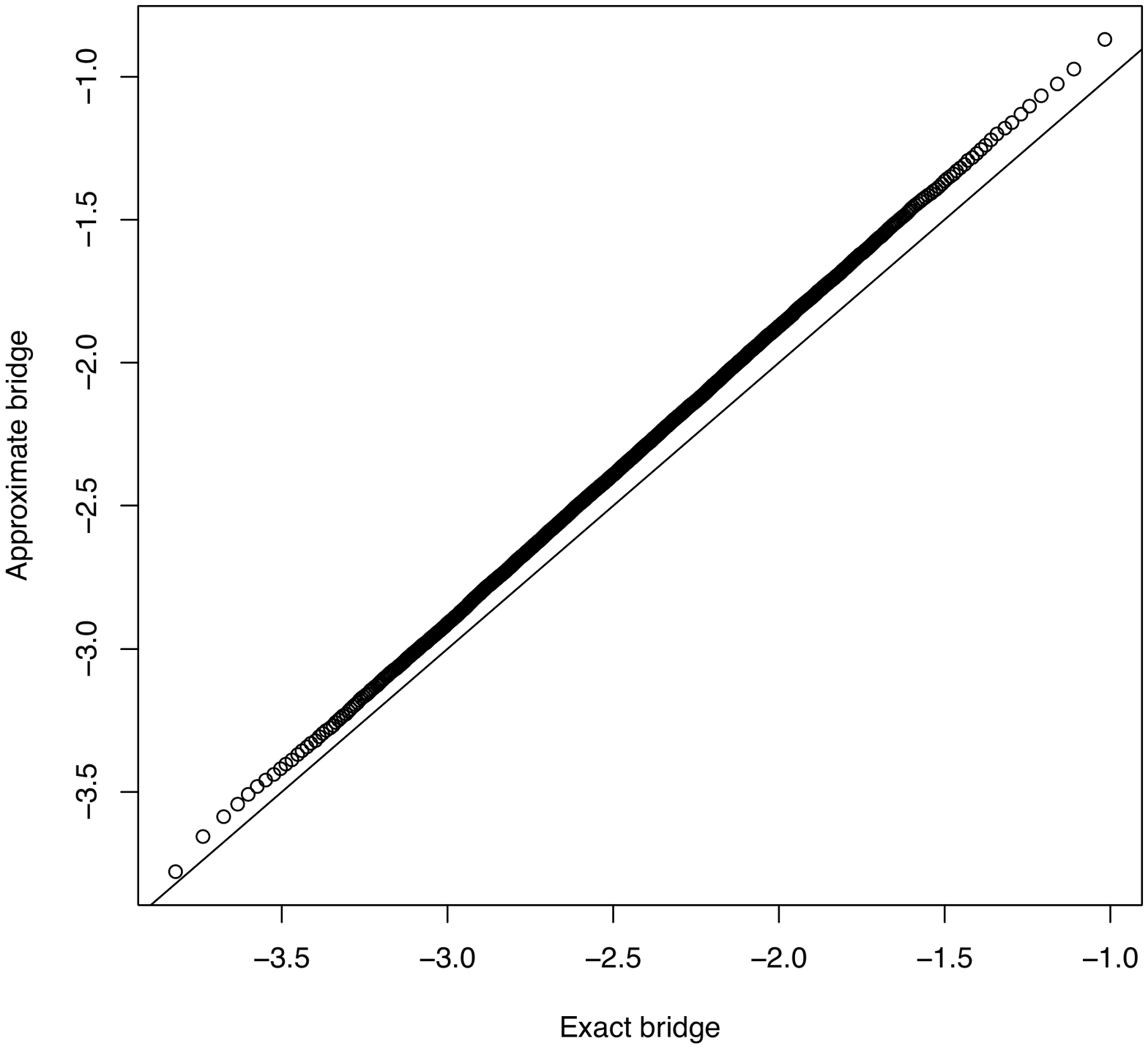}
\vspace{-5mm} 
\includegraphics[width=5.5cm]{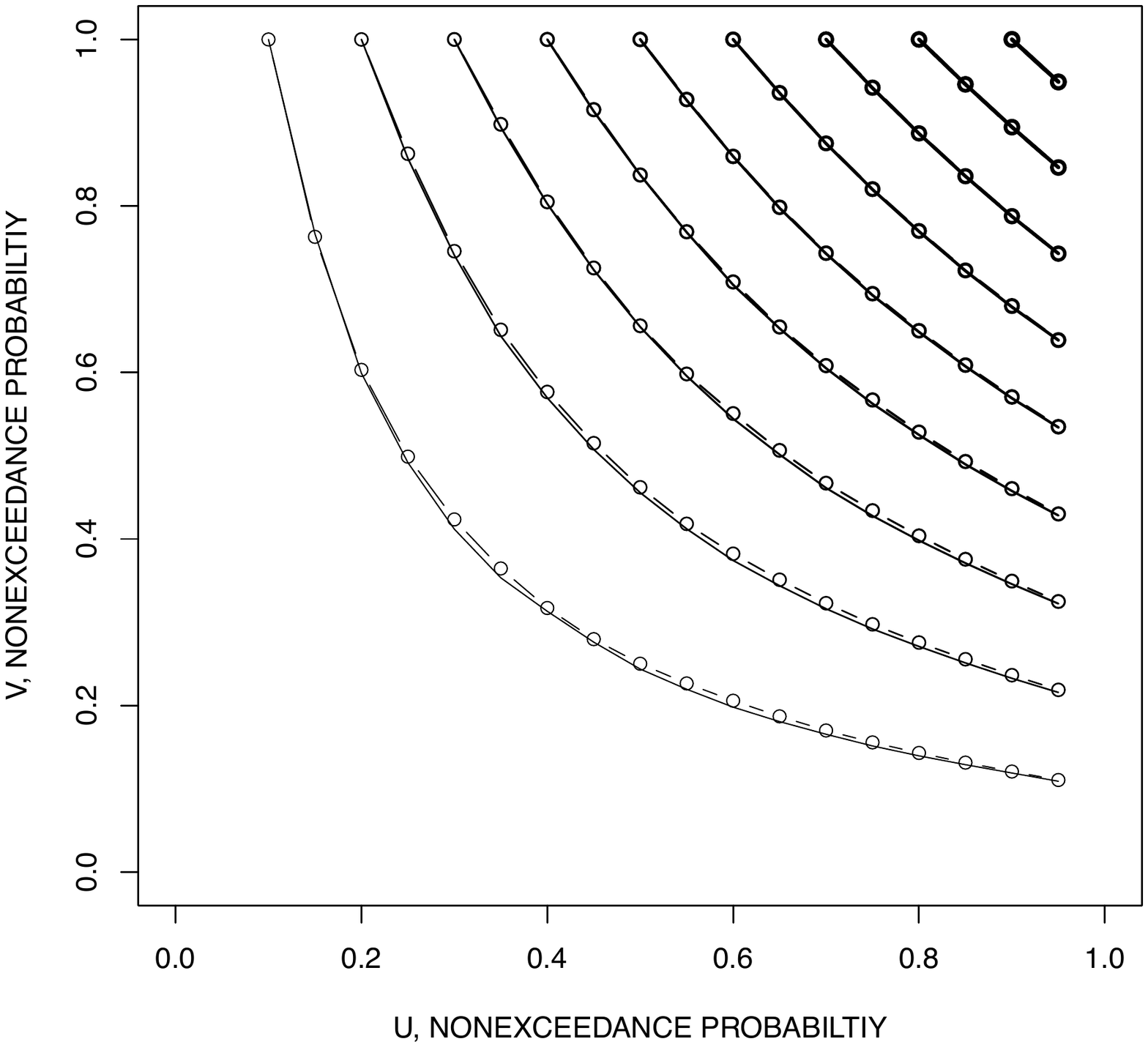} 
\vspace{-1.3cm}
\caption{\label{simulation7} 
Q-Q plots comparing the empirical marginal distributions at time 0.5 for 
50.000 two-dimensional Ornstein--Uhlenbeck bridges from $(2,-2)$ 
to $(3,-3)$ simulated by the approximate method with $\gamma = 0.5$ to
the exact marginal distributions of Ornstein--Uhlenbeck bridges. Level 
curves of the empirical copula for the two 2-dimensional distribution 
at time 0.5 are compared to 
to those of the exact copula (full drawn curves).} 
\end{center}
\end{figure}

\begin{figure}
\begin{center}
\vspace{-3.0cm}
\includegraphics[width=5.5cm]{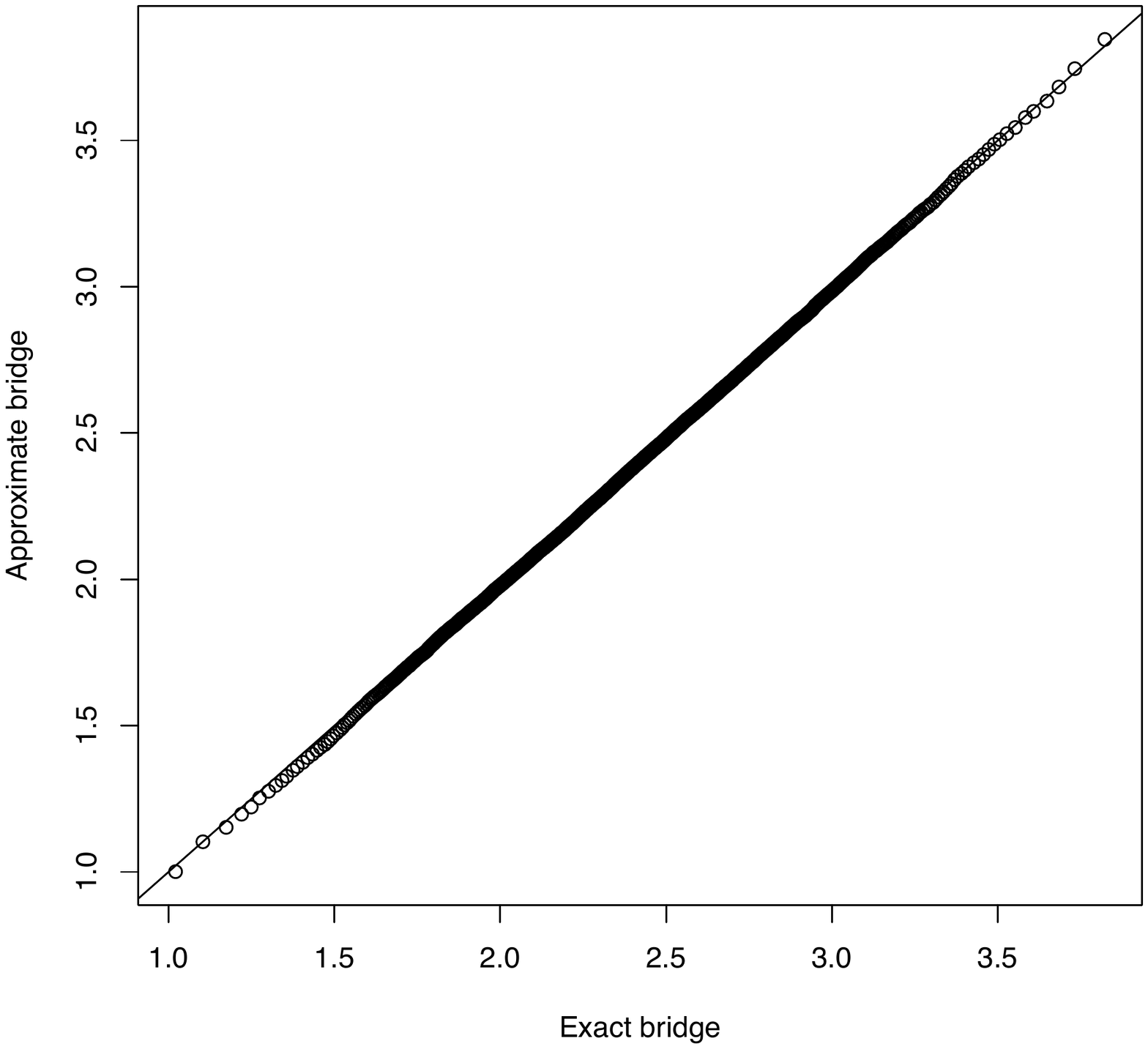} 
\hspace{-5mm}
\includegraphics[width=5.5cm]{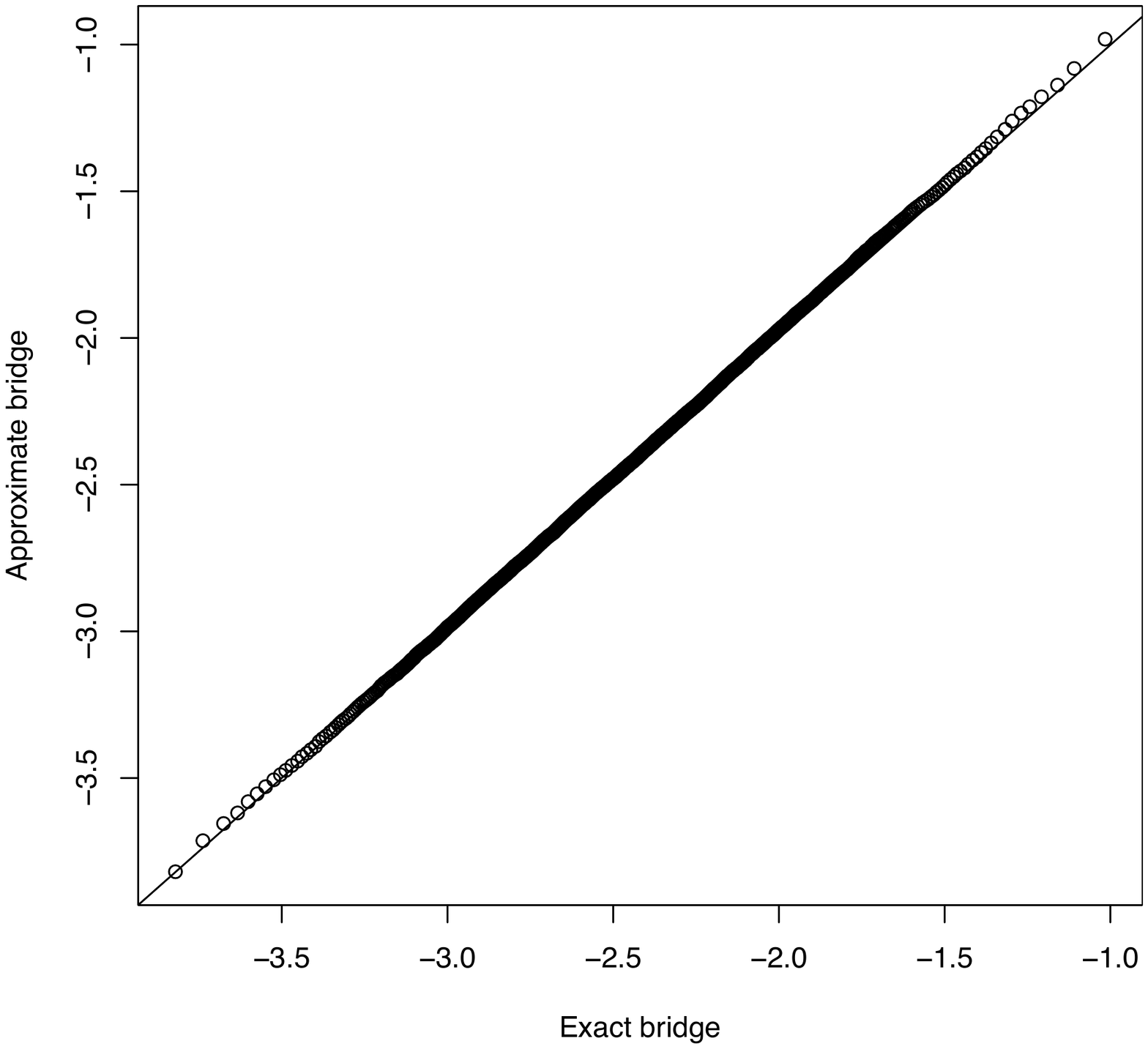}
\vspace{-5mm}
\includegraphics[width=5.5cm]{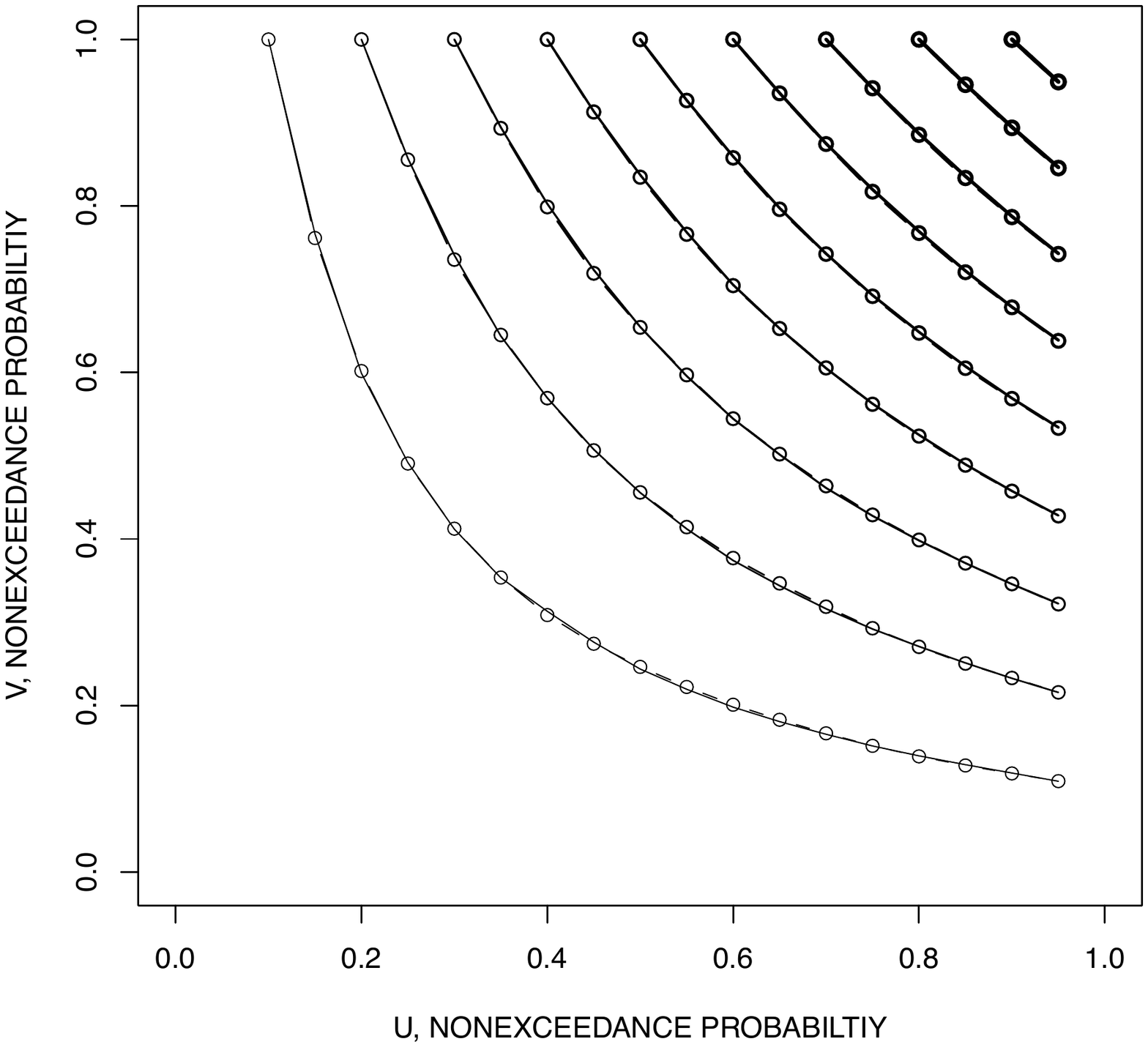} 
\vspace{-1.3cm}
\caption{\label{simulation8} 
Q-Q plots comparing the empirical marginal distributions at time 0.5 for 
50.000 two-dimensional Ornstein--Uhlenbeck bridges from $(2,-2)$ 
to $(3,-3)$ simulated by the approximate method with $\gamma = 0.9$ to
the exact marginal distributions of Ornstein--Uhlenbeck bridges. Level 
curves of the empirical copula for the two 2-dimensional distribution 
at time 0.5 are compared
to those of the exact copula (full drawn curves).} 
\end{center}
\end{figure}

Finally, diffusion bridges from $(2,-2)$ to $(3,-3)$ were simulated. These
are bridges from the boundary of the $98.2\%$-ellipse of the
stationary distribution to the boundary of the $99.99\%$-ellipse. 
Such a bridge is extremely rare in data.
We simulated 50.000 diffusion bridges using the approximate method
in Section \ref{approxbridge} with $\gamma = 0.5$, $\gamma = 0.9$ and
$\gamma = 0.99$. Computing times were in the interval 10 - 40  seconds
(increasing with the value of $\gamma$). 
Marginal distributions and the copula (at time 0.5) 
are compared to exact results (Lemma \ref{lemma:OU1}) in 
Figures \ref{simulation7} and \ref{simulation8}.
Even for this extreme bridge, the approximate diffusion bridge produces 
a surprisingly good fit to the distribution of the exact 
Ornstein--Uhlenbeck bridge for the coupling methods used here. The fit 
increases with $\gamma$. For $\gamma = 0.9$ the fit is 
almost perfect. For $\gamma = 0.99$ the fit is similarly
excellent. Both for this and the previous bridge it seems that  
the approximate method tends to become exact as $\gamma$ tends to one.
It is an intriguing question that requires further research, whether
this can be proved mathematically and in what generality it is true.

Here we restrict ourselves to a further simulation study in order to understand 
what happens when the approximate bridge simulation method gives almost
exact results for $\gamma$ close to one. From (\ref{z-density}) it
follows that if the function $\pi_T(x)$ is constant, then the
approximate method gives exact diffusion bridges. To investigate
whether $\pi_T(x)$ is in some cases constant (or varies only a
little), we ran (for each of a number of bridges and $\gamma$ values)
the pseudo-marginal MH-algorithm with 9000 iterations for $N = 50,
100, 150, 200$ and 300 (after a burn-in of 1000 iterations). For each
bridge, $\gamma$-value and $N$-value, the empirical mean and variance
were calculated of the values of $\hat \rho_N ({\bf T}^{(i)})$ from
all iterations in the run ($i = 1001, \ldots, 10000$). From results
on conditional expectations and the geometric distribution (and because
$\hat \rho_N$ is an unbiased estimator of $1/\pi_T(x)$ conditionally on
$X^{(i)} = x$)
\be
\label{regression}
V(\hat \rho_N) = V\left(\frac{1}{\pi_T(B)}\right) + E \left(
  \frac{1-\pi_T(B)}{\pi_T(B)^2} \right) \cdot \frac{1}{N},
\ee 
where $B$ in an $(a,b,T)$-bridge. Thus by linear regression of the
empirical variances of $\hat \rho_N$ from the MH-runs on $1/N$, we can
estimate $V(1/\pi_T(B))$. If the variance is zero, the function
$\pi_T(x)$ is constant, and in this case we can estimate the constant
value $\pi_T$ by the reciprocal empirical mean of the simulated $\hat
\rho$-values. The estimated slope can, in this case, be compared to
the value calculated by (\ref{regression}) using the estimated value
of $\pi_T$ as a consistency check. This check supported the
conclusions below.

For the unlikely bridge from $(0.785,0.785)$ to $(1.091,1.091)$ the
variance of $1/\pi_T(B)$ was not zero, but the ratio of the standard
deviation to the mean of $\hat \rho$ was 0.16 for $\gamma = 0.9$ and
0.04 for $\gamma = 0.99$, so for these $\gamma$-values $\pi_T(x)$ does
not vary much. For $\gamma=0$ and $\gamma = 0.5$ the standard
deviation was found to be of the same order as the mean.
For the very unlikely bridge from $(2,-2)$ to $(3,-3)$ the estimated
variance of $1/\pi_T(B)$ was zero for $\gamma = 0.9$ and $\gamma =
0.99$, implying an exact bridge. For $\gamma = 0.5$ the variance was
not zero, but the ratio of the standard deviation to the mean of $\hat
\rho$ was 0.06.

The simulation study indicates that the reason why exact or almost
exact diffusion bridges can be obtained by the approximate method when
$\gamma$ is close to one is that in this case $\pi_T(x)$ is constant
or almost constant. However, the simulation study also showed that the
probability $\pi_T(x)$ in such cases is very small, e.g.\ 0.03 for the
bridge from $(2,-2)$ to $(3,-3)$ for $\gamma = 0.99$.
For the Ornstein-Uhlenbeck process (\ref{OU}) we have
\[
d(X_t - X'_t) = -B(X_t - X'_t)dt + (1-\gamma)\sigma\Pi(X_t, X_t')dW_t 
- \sqrt{1-\gamma^2}\sigma u(X_t,X_t') dU_t. 
\]
With $\varepsilon = 1-\gamma$, we have $\sqrt{1-\gamma^2} \approx 
\sqrt{2\varepsilon}$ , so
\[
d(X_t - X'_t) \approx -B(X_t - X'_t)dt + \varepsilon\sigma\Pi(X_t,
X_t')dW_t - \sqrt{2\varepsilon}\sigma u(X_t,X_t') dU_t.
\]
For $\gamma$ close to one, $\varepsilon$ is small and $\sqrt{2\varepsilon} \gg \varepsilon$, so
the largest contribution is from the drift term, which should bring
$X_t$ and $X_t'$ close together, while
the contribution from $dW_t$ is very small compared to the contribution from $dU_t$.
Therefore if $X_t$ is an approximate $(a,b, T)$-bridge and $X'_t$ is
started from $A$, then the process $X_t - X'_t$
is nearly deterministic, but has a very small random contribution from
$U_t$ which is unlikely to bring the processes together, but is
independent of $X_t$. The contribution from $W_t$ is the only part
that depends on $X_t$, but is small even compared to the 
$U_t$ contribution, so it makes almost no difference to the path $X_t - X'_t$.
This seems a likely explanation why in this case $\pi_T(x)$ is small but almost constant. 
The probability of hitting the bridge depends mainly on the start point $A$ and a little on the process $U_t$
both of which are chosen independently of $X_t$. The contribution from the bridge itself, 
through $W_t$, is negligible.

The optimal choice of $\gamma$ is an interesting open question. The
solution might well depend on both the stochastic differential
equation and on the end-points of the bridge. For the approximate
method a reasonable definition of an optimal $\gamma$ value is the one
that provides the best fit to an exact bridge. For the exact MCMC
methods, the optimal $\gamma$-value is the one for which the computing
time is minimized. A possible solution is to draw $\gamma$ randomly in
the interval $[-1,1)$. In this case the question is what is the
optimal distribution from which to draw $\gamma$ (also model and
end-point dependent).

\sect{Bayesian estimation for discretely observed multivariate diffusions}
\label{bayes}

In this section we demonstrate how our method can be used for Bayesian
estimation for discretely observed multivariate diffusions. Specifically,
we consider estimation for the multivariate hyperbolic diffusion. 

The $d$-dimensional hyperbolic diffusion is given by
\be
\label{hyperbolicsde}
dX_t = -\frac{\alpha X_t}{\sqrt{1+\| X_t \|^2}} dt + dW_t,
\ee
where $\alpha>0$. It is the characteristic diffusion in the sense of
\cite{kent} for the multivariate hyperbolic distribution with density
function
\be
\label{hypdensity}
\nu(x) = \frac{(\alpha/\pi)^{\frac12 (d-1)}}{2K_{\frac12(d+1)}(2 \alpha)}
\exp \left( -2\alpha\sqrt{1+\| x \|^2} \right).
\ee
Here $K$ is a modified Bessel function of the third kind. The
hyperbolic distributions were introduced by \cite{oebn77}, who also
introduced the hyperbolic diffusions. The hyperbolic distribution has
heavier tails than the normal distribution. The
multivariate hyperbolic diffusion (\ref{hyperbolicsde}) is ergodic and
time-reversible with stationary density function $\nu(x)$; see
Section 10 in \cite{kent}.

Suppose we have observations at the time points $t_0=0 < t_1 < \cdots < 
t_n)$ from a hyperbolic diffusion. Then our data is a set of partial
observations of the full data set consisting of the continuous sample
path in the time interval $[0,t_n]$. We can therefore apply the Gibb's
sampler to generate draws from the posterior distribution. To do so we
need to be able to simulate the full continuous sample path
conditionally on the data $D = (X_{t_0}, \ldots , X_{t_n})$ and on $\alpha$. 
This is done by simulating independent hyperbolic diffusion bridges 
between the observations $X_{t_{i-1}}$ 
and $X_{t_{i}}$ in all intervals $[t_{i-1},t_i]$, $i=1,\ldots,n$.

We also need to simulate draws from the conditions distribution of
$\alpha$ given a continuous sample path. The likelihood function when
the data is a continuous sample path is given by Girsanov's
formula. The following expression without stochastic integrals for the
likelihood function can be obtained by applying Ito's formula to the
function $\sqrt{1+\| x \|^2} $, for details see Section \ref{proofs}.
\be
\label{likelihood}
L_{t_n}^c (\alpha) = 
\exp \left( \alpha H_{t_n} - \frac12 \alpha^2 B_{t_n} \right),
\ee
where
\[
H_t = \sqrt{1+\| X_0 \|^2} - \sqrt{1+\| X_t \|^2} +\int_0^t 
\frac{1+\frac12 \| X_s\|^2}{(1+\| X_s\|^2)^{3/2}} ds 
\]
and 
\[
B_t = \int_0^t \frac{\| X_s\|^2}{1+\| X_s\|^2} ds 
\]

The continuous time model is an exponential family of stochastic processes 
in the sense of \cite{kuso}. The conjugate prior is a normal distribution. 
If we choose as our prior the normal distribution with expectation 
$\bar \alpha$ and variance $\sigma^2$, the posterior distribution is a 
normal distribution with expectation 
$(H_{t_n}+\bar \alpha/\sigma^2)/(B_{t_n}+\sigma^{-2})$ and variance 
$(B_{t_n} + \sigma^{-2})^{-1}$.

\subsection{Simulations}

\begin{figure}
\begin{center}
\vspace{-1cm}
\includegraphics[width=6cm]{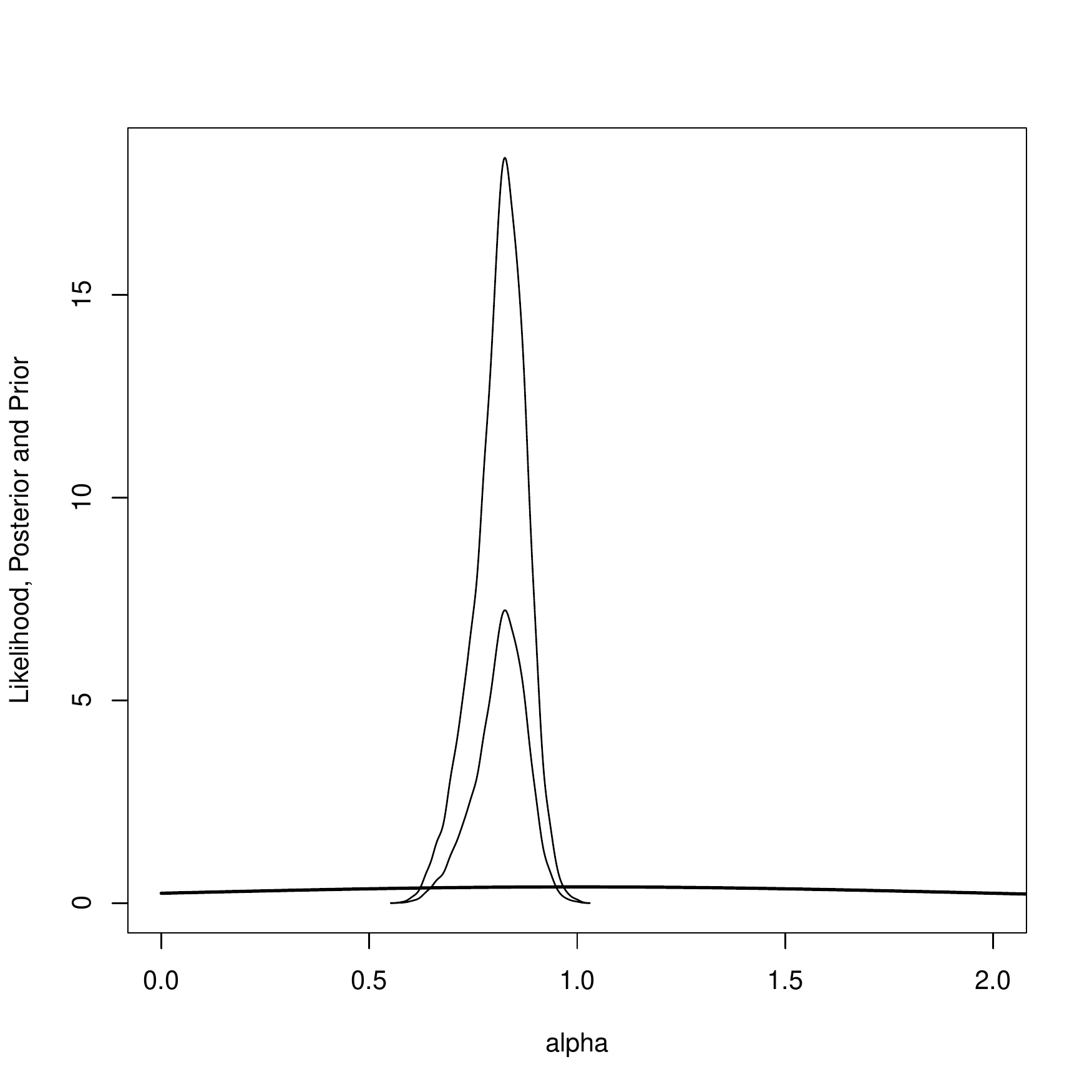} 
\hspace{15mm}
\includegraphics[width=6cm]{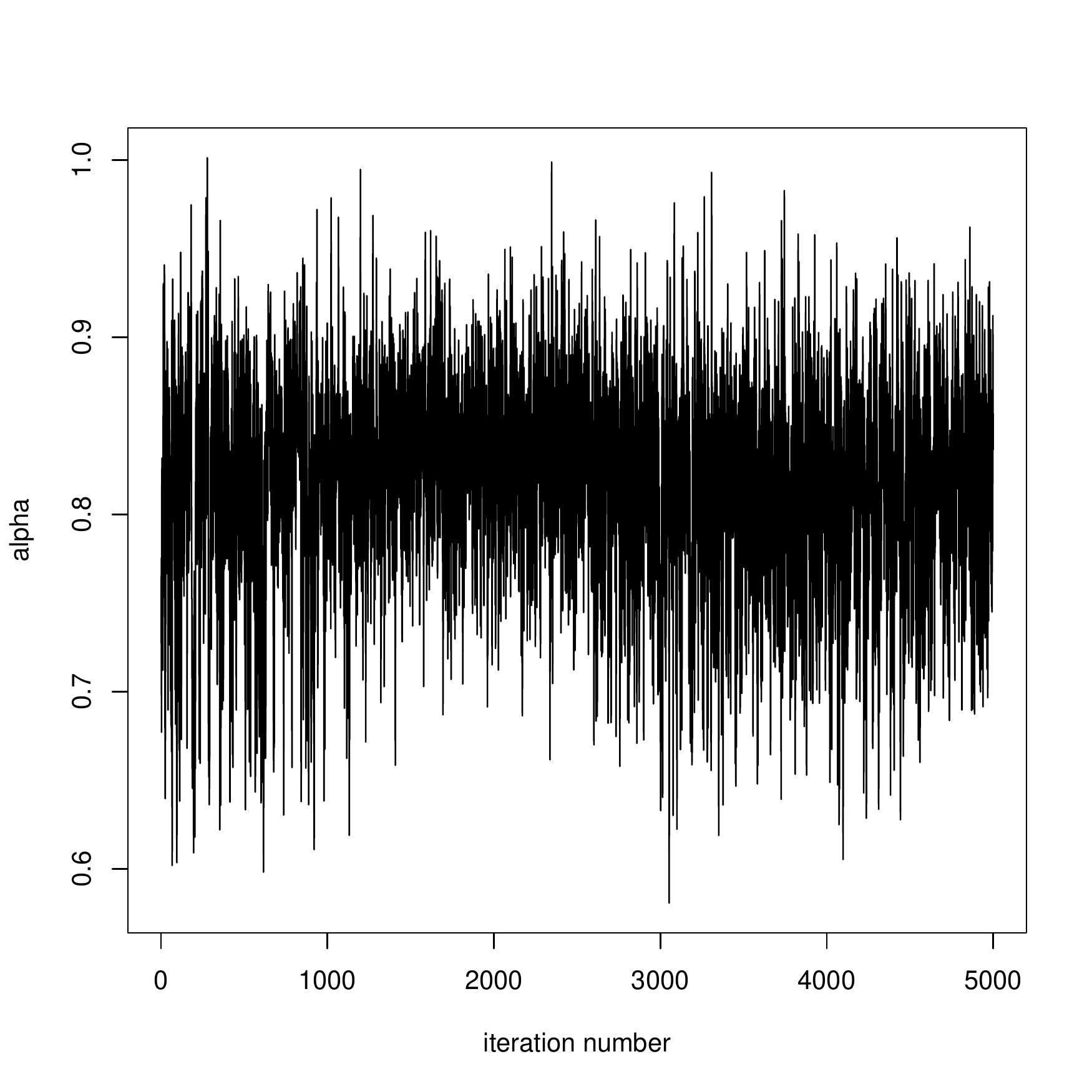}  \\
\caption{\label{estimation} 
The plot to the left is the posterior, the prior and the likelihood 
function for a sample of 1000
observations of the two-dimensional hyperbolic diffusion with $\alpha
= 0.8$ obtained by 5000 iterations of the Gibbs sampler using the the
approximate method in Section \ref{approxbridge}  with $\gamma =
0.5$. To the right is shown the time series of the 5000 draws of $\alpha$.} 
\end{center}
\end{figure}

In order to test how well our bridge simulation method works for
estimation of parameters in multivariate diffusions, we simulated a
sample of observations at the time points $t_i = i$,
$i=1,\ldots, 1000$ of the two-dimensional hyperbolic diffusion with
$\alpha = 0.8$. As prior we used the $N(1,1)$-distribution.
Then we ran 5000 iterations of the Gibbs sampler that starts by drawing
$\alpha$ from the prior and then alternates between drawing a
continuous sample path in $[0,1000]$ given the data and $\alpha$ and
drawing $\alpha$ from the conditional distribution given the
continuous sample path. To simulate the continuous sample path given
the data, we used the approximate simulation method in Section
\ref{approxbridge} with $\gamma = 0.5$. Figure \ref{estimation} shows
the prior distribution, the posterior distribution, and the likelihood 
function. It also shows a plot of the time series of draws of
$\alpha$ for the 5000 iterations that indicates that the algorithm 
stabilizes very quickly. The posterior distribution is 
nicely concentrated around the true parameter value $\alpha = 0.8$. 
The 95\%-credibility interval is $[0.6783, 0.9239]$ and the mean 
posterior estimate is 0.8164.

\sect{Proofs}
\label{proofs}

\noindent
{\bf Proof of Lemma \ref{lemma:samplepaths}:} The results
(\ref{prime-dif}) is straightforward. In order to prove the result on
the sample path of $X$,  we begin by introducing
a new Wiener process to make the structure of the Wiener process $W'$ clearer.
First choose (for each pair $(x,x')$) an orthonormal base $u_2(x,x'),
\ldots u_d(x,x')$ for the orthogonal complement of the space spanned 
by $u(x,x')$. This can obviously be done such that $u_i(x,x')$ is a
continuous function of $(x,x')$. Define
\[
V^i_t = \int_0^t u_i(X_s,X'_s)^T dW'_s, \ \ i=1,\ldots,d,
\]
where $u_1=u$. Clearly, $V_t = (V^1_t,\ldots,V^d_t)$ is a
$d$-dimensional standard Wiener process, and $dW'_t = \sum_{i=1}^d 
u_i(X_s,X'_s) dV^i_t$. Since $\Pi(X_t,X'_t)$ is the 
projection on the space spanned by $u(X_t,X'_t)$, while $I-\Pi(X_t,X'_t)$ 
is the projection on the orthogonal complement to this space, it follows that
\[
V^i_t = \int_0^t u_i(X_s,X'_s)^T{\cal O}(X_s,X'_s) dW_s, \ \ i=2,\ldots,d,
\]
while
\[
V^1_t = \gamma \tilde Y_t + \sqrt{1-\gamma^2} U_t,
\]
where
\[
\tilde Y_t = \int_0^t u(X_s,X'_s)^T {\cal O}(X_s,X'_s) dW_s
\]
is a standard Wiener process. The process
\[
Y_t = \tilde Y_t - \gamma V^1_t
\]
is a Wiener process with infinitesimal variance $1-\gamma^2$. It is 
independent of the Wiener process $V_t$, because its 
quadratic co-characteristics with the components of $V$ are all zero. For 
instance, $\langle Y, V^1\rangle_t = \gamma (\langle \tilde Y 
\rangle_t - \langle V^1\rangle_t) = 0$. It follows that $Y$ is 
also independent of the Wiener process $W'$ that drives $X'$, again 
because the quadratic co-characteristics are zero. When $\gamma = -1$,
$Y_t =0$, and when $\gamma = 0$, $Y_t = \tilde Y_T$ and $V^1_t = U_t$.

Since $[ I-\Pi(X_t,X'_t)]dW'_t = [ I-\Pi(X_t,X'_t)] {\cal O}(X_s,X'_s) dW_s$ 
and $\tilde Y_t = Y_t + \gamma V^1_t$, we find that
\bean
{\cal O}(X_s,X'_s) dW_t &=& \{ I-\Pi(X_t,X'_t)\}{\cal O}(X_s,X'_s) dW_s 
+ \Pi(X_t,X'_t) {\cal O}(X_s,X'_s) dW_s \\
&=& \{ I-\Pi(X_t,X'_t)\}dW'_t + u(X_t,X'_t)d\tilde Y_t \\
&=& \{ I-(1-\gamma)\Pi(X_t,X'_t)\}dW'_t + u(X_t,X'_t)dY_t,
\eean
from which (\ref{stjerne-dif}) follows.

\halmos

\noindent
{\bf Proof of Theorem \ref{theorem1}:} By the strong Markov property
$Z$ has the same distribution as $X'$, so the conditional distribution
of $\{ Z_t \}_{0\leq t \leq T}$ given $Z_T = b$ is the distribution of a 
$(a,b,T)$-diffusion bridge. Now 
\[
P(Z \in \cdot \, | \, X_T = b, \tau \leq T) = P(Z \in \cdot \, | \,
Z_T = b, \tau \leq T),
\]
and the event $\{ Z_ = b, \tau \leq T \}$ is the event that $Z$ is a
$(a,b,T)$-diffusion bridge and that the diffusion bridge is hit by
$X$. The theorem follows because the sample path of $X$ up to time $\tau$ 
is $\tilde {\cal K}_\tau \left( X_0 , \{ Z_s \}_{0 \leq s \leq \tau}, 
\{ U'_s \}_{0 \leq s  \leq \tau} \right)$. The distribution of $X_0$ 
conditional on $X_T=b$) has density $p^*_T(b, \cdot)$. A proof of the
last claim is given  in the proof of Lemma \ref{lemma:fundamental}. 

\halmos

\noindent
{\bf Proof of Lemma \ref{lemma:fundamental}:} The second identity in 
(\ref{bridgetransition}) follows from (\ref{balance}). The first expression 
for $q$ is the well-known expression for the transition density of a
diffusion bridge ending in $b$ at time T, see \cite{yor}, p.\ 111. It
can be easily established by direct calculation. 
The second expression for $q$ can similarly be obtained
as the transition density of $\bar X$ by direct calculation: the
conditional density of $\bar X_t$ given $\bar X_s$ ($s<t$) is
\[
p_{\bar X_s, \bar X_t}(x,y)/p_{\bar X_s}(x) = 
p_{X^*_{T-t},X^*_{T-s}}(y,x)/p_{X^*_{T-s}}(x) = p^*_{T-t}(b,y)
p^*_{t-s}(y,x)/p^*_{T-s}(b,x).
\]

Now assume that $X_0 \sim \nu$. Then $X_T \sim \nu$, and the joint
density of $(X_0, X_T)$ is $\nu (x) p_T(x,y)$ $ = \nu(y) p^*_T(y,x)$, 
again by (\ref{balance}). Hence the conditional density of $X_0$ given
$X_T=b$ is $p^*_T(b,x)$. Obviously, the density of $\bar X_0 = X^*_T$ 
is $p^*_T(b,x)$, so the process $\{ \bar X_t  \}$ and the conditional 
process $\{ X_t \}$ given that $X_T=b$ (both of which are Markov
processes) have the same transition densities and the same initial
distribution. Therefore they have the same distribution.

\halmos

\noindent
{\bf Proof of Corollary \ref{cor3}:} The function $\tilde {\cal
  K}_T \left( A, x, U \right)$ 
is not defined for all $x \in C_T$, but it is defined for all relevant
trajectories $x$. For other $x$ it can be given an arbitrary definition.
To prove equation (\ref{z-density}) note that the joint density of a
diffusion bridge $X$, the independent random variable $A$, and an 
independent Wiener process $U$ conditional on the event that $X$ 
intersects with $\tilde {\cal K}_T \left( A, X, U \right)$ is 
\[
f_{br}(x) p^*_T(b,a) f_W(u) 1_M(x,\tilde {\cal K}_T \left( a, x, u 
\right))/\pi_T.
\]
Here $f_W$ is the density on $C_T$ of a standard Wiener process
w.r.t.\ a suitable dominating measure. From this expression
(\ref{z-density}) follows by marginalization.

\halmos

\noindent
{\bf Proof of Lemma \ref{lemma:michael}:}
The result follows straightforwardly because $(Z_{t_1},
\ldots , Z_{t_n})$ is a linear transformation of a multivariate normal
distribution. For completeness, we give the details. 
The distribution of $(X_{t_1}, \ldots , X_{t_{n+1}})$ equals the
distribution of observations from an Ornstein--Uhlenbeck process
started at $X_0=x_0$. Define
\[
Y_{t_i} = X_{t_i} - e^{-B(t_{n+1} - t_{i})} \Gamma_{t_i} 
\Gamma_{t_{n+1}}^{-1} X_{t_{n+1}}, \ \ i=1, \ldots, n.
\]
It is well-known (and easily checked) that $(Y_{t_1}, \ldots ,
Y_{t_{n}})$ is independent of $X_{t_{n+1}}$. Since $Z_{t_i} = Y_{t_i}
+ e^{-B(t_{n+1} - t_{i})} \Gamma_{t_i} \Gamma_{t_{n+1}}^{-1} x$, $i=1,
\ldots, n$, it follows that the distribution of $(Z_{t_1}, \ldots ,
Z_{t_{n}})$ equals the conditional distribution of $(X_{t_1}, \ldots ,
X_{t_{n}})$ given $X_{t_{n+1}}=x$.

\halmos

\noindent
{\bf Proof of (\ref{likelihood}):}

By Girsanov's formula, the likelihood function is given by
(\ref{likelihood}) with $B_T$ as in Section \ref{bayes} and 
\[
H_T = - \int_0^T \frac{X_s^T}{\sqrt{1+\| X_s \|^2}} dX_s,
\]
see e.g.\ \cite{kuso}, p.\ 297. By applying Ito's formula to the function $F(x) 
= \sqrt{1+\| x \|^2} $, $x \in \R^2$, we see that
\[
F(X_T) = F(X_0) + \int_0^T \frac{X_s^T}{\sqrt{1+\| X_s \|^2}} dX_s
+ \int_0^T \frac{1+\frac12 \| X_s\|^2}{(1+\| X_s\|^2)^{3/2}} ds.
\]
From this the expression for $H_T$ in  Section \ref{bayes} follows.

\section*{Acknowledgement} 

Mogens Bladt acknowledges the support from grant SNI15945 by the
Mexican Research Council CONACYT.  The research of Michael S\o rensen
was supported by the Center for Research in Econometric Analysis of
Time Series funded by the Danish National Research Foundation, and by
two grants from the University of Copenhagen Programme of Excellence.

\bibliographystyle{natbib}

\end{document}